\documentclass[right]{amsart}
\usepackage{myarticle}

\usepackage{graphicx}
\graphicspath{{SpaceTimePlots/}}
\usepackage{amsbsy}
\usepackage{latexsym}
\usepackage{amssymb}
\usepackage{amsmath}
\usepackage{float}
\usepackage[latin1]{inputenc}
\usepackage{textcomp}
\usepackage{pgfplots}
\usepackage{array}
\usepackage{color}
\usepackage{caption}
\usepackage{mathtools}
\usepackage{algorithm}
\usepackage[noend]{algorithmic}

\definecolor{lightgreen}{rgb}{0.6,1,1}

\usepackage{tikz}
\usetikzlibrary{shapes.geometric, arrows}

\tikzstyle{startstop} = [rectangle, rounded corners, minimum width=1cm, minimum height=0.5cm,text centered, draw=black, fill=red!30]
\tikzstyle{process} = [rectangle, minimum width=1cm, minimum height=0.5cm, text centered, draw=black, fill=orange!30]
\tikzstyle{decision} = [diamond, minimum width=1cm, minimum height=1cm, text centered, draw=black, fill=green!30]
\tikzstyle{io} = [trapezium, trapezium left angle=70, trapezium right angle=110, minimum width=3cm, minimum height=1cm, text centered, draw=black, fill=blue!30]
\tikzstyle{arrow} = [thick,->,>=stealth]

\def\vec#1{\ensuremath{\mathchoice
           {\mbox{\boldmath$\displaystyle\mathbf{#1}$}}
           {\mbox{\boldmath$\textstyle\mathbf{#1}$}}
           {\mbox{\boldmath$\scriptstyle\mathbf{#1}$}}
           {\mbox{\boldmath$\scriptscriptstyle\mathbf{#1}$}}}}

\newcommand{\bdot}[1]{{\boldsymbol{\dot{#1}}}}

\newcommand\id{\operatorname{id}}

\renewcommand\div{\operatorname*{div}}

\newcommand\trace{\operatorname*{trace}}

\newcommand\Real{{\mathbb R}}

\renewcommand\bar{\overline}%

\definecolor{orange}{RGB}{255,127,0}
\definecolor{dblue}{RGB}{16,78,139}
\definecolor{green}{rgb}{0.1, 0.6, 0.3}

\newcommand{\gruen}[1]{\textcolor{black}{#1}}

\def\tensor{\,\raise2pt\hbox{${}_{\otimes}$}\,}
\def\Paragraph#1{\smallskip\paragraph{\bf #1}\smallskip}

\usepackage[colorlinks=true,
            pdftitle={Computational dynamics for linear solids with fractures},
            pdfauthor={Kerstin Weinberg, Christian Wieners},
            pdfpagemode={UseOutlines},pdfstartview={Fit},
            linkcolor=blue,urlcolor=blue,citecolor=blue,
            pdfsubject={Computational dynamics for linear solids with fractures}]{hyperref}

\def\Clearpage{}

\begin{document}

\title{Dynamic phase-field fracture with a first-order\\
  discontinuous Galerkin method for elastic waves}

\author{Kerstin Weinberg}

\address{Chair of Solid Mechanics, University of Siegen, Siegen, Germany.
  Email: {\tt kerstin.weinberg@uni-siegen.de} }

\author{Christian Wieners}

\address{Institute  of Applied and Numerical Mathematics, KIT, Karlsruhe, Germany.
  Email: {\tt christian.wieners@kit.edu} }

\begin{abstract}
We present a new numerical approach for wave induced dynamic
fracture. The method is based on a discontinuous Galerkin
approximation of the first-order hyperbolic system for elastic waves and a
phase-field approximation of brittle fracture driven by the maximum
tension.  The algorithm is staggered in time and combines an implicit
midpoint rule for the wave propagation followed by an implicit Euler
step for the phase-field evolution. At fracture, the material is
degraded, and the waves are reflected at the diffusive interfaces.
Two and three-dimensional examples demonstrate the advantages of the
proposed method for the computation of crack growth and spalling
initiated by reflected and superposed waves.
\end{abstract}



\maketitle

\section{Introduction} \label{sec:intro}
%
%
Dynamic fracture and fragmentation of solids attract constant
attention from scientists, both out of simple curiosity, technological
interest, and numerical challenges. Consequently, numerous
computational fracture schemes have been designed, including classical
discontinuous methods, such as interface 
formulations or local enrichment strategies, and continuous methods
with discontinuities replaced by sharp gradients. The phase-field
approach to fracture, chosen here, falls into the latter category and
implicitly represents the cracks by an additional continuous field, an order
parameter indicating the transition from elasticity to fracture.
As a result, the fracture evolves naturally within the numerical computation.

Numerical schemes for dynamic phase-field fracture are well
established, see, e.g.,
\cite{Miehe2010,borden2012phase,hesch2014thermodynamically}, and, more recently,
\cite{ren2019explicit,mandal2020evaluation}, and the references therein
for an overview. They all have in common that they use the classical
description of motion derived from the conservation of linear
momentum, which results in a second-order partial differential
equation. Their integration in time is performed with finite
differences or a Newmark family method, which inevitably introduces
some numerical dissipation \cite{HughesB}. In sudden rupture or for
impact problems, this spurious energy loss is usually manageable or
insignificant. However, for the propagation of waves, their
reflection, transmission, and superposition, the transported energy is
significant, and a spurious dissipation no longer allows for a
meaningful evaluation of the results. Special integration techniques
try to circumvent these spurious effects, but ultimately energy loss
is always a problem.

In order to solve this problem, we propose here an approach that has
already proven successful in calculating conservation laws. For this
purpose, the elastic equation of motion is transformed into a coupled
system of first-order differential equations, suitably discretized,
and combined with a phase-field method to calculate the
crack initiation and propagation.
\gruen{The use of first-order discontinuous Galerkin methods combined with implicit time-stepping methods for the elastic wave propagation is motivated by the high accuracy, the energy conservation as long no cracks open, and the small numerical dissipation of this method. We show in our examples that these properties are essential to simulate fracture initiated by the superposition and reflection of waves. Nevertheless, these effects can be computed only with a high numerical expense on fine meshes, so that efficient parallel finite element software, including fast iterative solution methods in every implicit time step, is an indispensable requirement.
}

This paper is organized as follows. In
Section~\ref{sect:phasefieldfracture}, we define a phase-field model
for dynamic fracture, where the phase-field evolution is driven by a
stress-based fracture criterion \cite{bilgen2019crack}
and where the crack is approximated by a material degradation
depending on the phase field. For the discrete phase-field evolution,
we use lowest order conforming finite elements in space and the
implicit Euler method in time. Then, in Section~\ref{sec:DG}, we
recall the discontinuous Galerkin (DG) approximation for linear waves
\cite{HesthavenWarburton08}, based on the formulation of a first-order
hyperbolic system and using full upwind flux \cite{HocPSTW15}.  This
DG approximation is combined with the energy-conserving implicit
midpoint rule in time. Then, in Section~\ref{sect:phasefieldfracture},
we introduce a fully coupled algorithm combining the DG
method with the phase field and fracture evolution in a staggered
scheme with variable time steps depending on the fracture evolution.


\gruen{
Our numerical experiments in Section~\ref{sec:numerics} are motivated
by an Hopkinson bar experiment. In the first example}, we
study a configuration in two dimensions, where the crack is initiated
by the reflection and superposition of two incoming pressure waves. We
show that our staggered scheme
reproduces identical dynamic crack
behavior on different meshes, indicating convergence of the numerical
method. Then, this test is transferred to a three-dimensional
configuration which also shows crack opening by the superposition of
two incoming wave signals. All numerical results are realized within
the parallel finite element system~M++
\cite{Wieners10,baumgarten2020parallel}.  In
Section~\ref{sec:fracure}, we conclude with a discussion on
variational formulations of the phase-field driving forces and
possible extensions of our computational approach 
to dynamic models for non-linear elastic and dispersive waves.


\section{A phase-field approach for dynamic fracture}\label{sect:phasefieldfracture}

Physically, the initiation and propagation of cracks depend on various
conditions, such as the geometry of the structure, the material
toughness, the loading rate, the loading magnitude and its
distribution. In classical fracture mechanics, the stress state at the
crack tip determines the steady or unsteady crack growth
\cite{AndersonB}. Following the brittle fracture approaches of
Griffith and Irwin \cite{Griffith1921,Irwin1958}, the material fails
when the energy release reaches a critical value $G_\text{c}$.  This
assumption provides a criterion for crack growth; however, it needs to
be embedded into an energy minimization setting for the entire
structure to determine crack propagation, paths, and branching.

The non-local approach using a phase field goes back to the
fundamental variational model of Francfort and Marigo
\cite{FrancfortMarigo1998,bourdin2008variational} and has gained much
attention recently,
cf.~\cite{HenryLevine2004,Miehe2010,borden2012phase,deLorenzi2016,Negri2016}. The
phase-field crack modeling is popular mainly because of its potential
for capturing the evolution of complex crack patterns without the need
for specialized crack-front tracking algorithms.

In this section, we introduce the notation, and we summarize the model
and its phase-field approximation, which builds the basis for our
numerical realization.

\subsection{A dynamic model for crack evolution at small strains}

We consider a bounded Lipschitz domain $\Omega\subset \Real^d$ with
boundary $\partial \Omega= \bar{\Omega}\setminus\Omega$
and the time interval $[0,T]$.
We want to determine the evolution of a crack in time
$t\longmapsto \Gamma_\text{c}(t)\subset \Omega$ together with the dynamics of
elastic waves described by the velocity
vector $\vec v$ and the stress tensor $\vec \sigma$
      \begin{align}
      \label{eq:elastic}
      \varrho \bdot {\vec v}(t) = \div   \vec \sigma(t) + \vec f(t)
      \,,\qquad
      \bdot {\vec\sigma}(t) =
      \vec C\vec \varepsilon\big(\vec v(t)\big)
      \qquad\text{in }\Omega \setminus\Gamma_\text{c}(t) \,.
    \end{align}
Here, $\varrho$ is the mass density, $\vec f(t)$ the applied
volume force and $\vec C$ the Hookean elasticity tensor with
Lam\'{e} constants $\lambda$ and $\mu$.
Furthermore, for the combined model we assume that
\begin{itemize}
  \item
    the crack evolution is irreversible, i.e.,
    $\Gamma_\text{c}(t_1) \subset \Gamma_\text{c}(t_2)\subset\Omega$ for $t_1<t_2$;
  \item
    the elastic solution remains admissible with respect to a
    suitable fracture criterion;
  \item
    the crack evolution is determined by an energy dissipation criterion.
\end{itemize}

The wave propagation described by \eqref{eq:elastic} is driven by
initial and boundary conditions on $\partial\Omega =
\partial_\text{N}\Omega \cup \partial_\text{D}\Omega$
together with  free Neumann boundary conditions at the crack interface, i.e.,
      \begin{subequations}
      \begin{alignat}{2}
      \label{eq:elastic_2}
            {\vec v}(0) &= \vec v_0
            &&\ \ \
            \text{ and }\ \
                   {\vec \sigma}(0) = \vec \sigma_0
                   \ \text{ in }\Omega\,,\\
            \vec\sigma(t) \vec n &=\vec g_\text{N}(t)
            \text{ on }\partial_\text{N}\Omega
            &&\ \ \
            \text{ and }\ \
            \vec v(t)  = \vec v_\text{D}(t)
            \text{ on }\partial_\text{D}\Omega
            \,,
            \\
      \vec\sigma(t) \vec n &=\vec 0
      \ \text{on }\Gamma_\text{c}(t)\,.
    \end{alignat}
      \end{subequations}
The configuration depends on initial data $\vec v_0$ and $\vec \sigma_0$,
volume forces $\vec f$ in $\Omega$, and
boundary data $\vec g_\text{N}$ and $\vec v_\text{D}$.

The model has to be complemented by an evolution law for $\Gamma_\text{c}$.

\clearpage

\Paragraph{Energy dissipation}
At time $t$, the \gruen{mechanical energy of the} first-order system is given by
    \begin{align}
      \label{eq:ee}
      \mathcal E(t) =
      \frac12
       \int_{\Omega\setminus\Gamma_\text{c}(t)}
       \Big(\varrho|\vec v(t)|^2
       +\vec\sigma(t)\cdot\vec C^{-1}\vec \sigma(t)\Big)\,\mathrm d\vec x
       \,.
    \end{align}
In the case of homogeneous data, $\vec f = \vec 0$ and $\vec g_\text{N} = \vec 0$,
we have $\mathcal E(t_1) = \mathcal E(t_2)$ if
$\Gamma_\text{c}(t_1) = \Gamma_\text{c}(t_2)$.
In the case of crack propagation, $\Gamma_\text{c}(t_2) \setminus \Gamma_\text{c}(t_1)\neq \emptyset$, we claim for
energy dissipation,
     \begin{align*}
       \mathcal D(t_1,t_2)
       &=
       \mathcal E(t_1) -  \mathcal E(t_2)
       \\
       &= \
       G_\text{c} \, |\Gamma_\text{c}(t_2)\setminus\Gamma_\text{c}(t_1)|_{d-1}
       \,,\qquad t_1<t_2\,,
     \end{align*}
i.e., the energy release is proportional to the (scaled) $d-1$ surface
volume of the crack increment, where the Griffiths constant
$G_\text{c} > 0 $ is a material parameter determining the critical energy release rate.
Energy dissipation in a dynamic phase-field
approximation can be included by a minimization formulation in every
time step \cite{larsen2010existence}. Here
we propose a different approach where energy dissipation is integrated in the time-stepping
scheme; details will be outlined in Section~\ref{sec:fracure}.


\Paragraph{Crack criteria}
To formulate a fracture evolution criterion we use here 
the \emph{maximum principal stress}
     \begin{align*}
       \sigma_\text{I} = \max_{\vec n\in S^2} \vec n\cdot \vec\sigma\vec n\,,\qquad
       S^2 = \big\{ \vec n\in \Real^3\colon |\vec n| = 1\big\}
     \end{align*}
and relate it to a material resistance value $\sigma_\text{c} > 0$.
We assume that the material response is elastic, if for all $\vec n\in S^2$ the tension
$\vec\sigma\vec n$ in direction $\vec n$ is smaller than
$\sigma_\text{c}$. Otherwise, if the maximum principal stress gets
larger, i.e., $\sigma_\text{I}(t,\vec x)\geq \sigma_\text{c}$, a driving force
is activated which initiates fracture.


\gruen{
Note that the stress-based criterion only
activates fracture driving forces in our dynamic model, and sufficiently strong forces are required to open a crack. This criterion is different from a quasi-static model with pointwise constraint $\sigma_\text{I}\leq \sigma_\text{c}$ which results in an instantaneous nonlinear material response and thus in a faster material response than wave propagation.
}

\subsection{A phase-field approximation}

The crack evolution is described by an order parameter, the phase field
     \begin{align*}
       s\colon [0,T] \times \overline \Omega \longrightarrow [0,1]\,,
     \end{align*}
starting with $s(0,\vec x) = 1$ for $\vec x\in \Omega$.  In our
approach, the phase-field evolution is reversible,
\gruen{
so that in the case that the
elastic driving force is not strong enough to initiate fracture,
the phase field recovers if the elastic driving force is getting smaller.}
The irreversible physical process of crack opening depends on the history, the
infimum phase field
    \begin{align*}
      s_\text{inf}(t,\vec x) &= \inf\big\{s(t',\vec x) \colon 0\leq t' \leq t\big\}
      \,,
    \end{align*}
determing the elastic domain
    \begin{align*}
      \Omega\left(s_\text{inf}\right)
      &= \big\{\vec x\in \Omega \colon
      s_\text{inf}(t,\vec x) > s_\text{min}\big\}
      \,.
    \end{align*}
Depending on the critical parameter $s_\text{min}> 0$
this defines the phase-field approximation
of the fracture zone $\Omega\setminus\Omega\left(s_\text{inf}\right)$.
Note that this parameter
is a regularization which extends the fracture $\Gamma_\text{c}(t)$ to a volume
$\Omega\setminus\Omega\left(s_\text{inf}\right)$.

Depending on the phase field $s$ at time $t\in [0,T]$,
the crack surface energy is approximated by
     \begin{align}
       \mathcal G(s) = \frac{G_\text{c}}{l_\text{c}}
       \int_\Omega \gamma_\text{c}(s)\, \mathrm d\vec x
       \,\ \text{with }\quad
       \gamma_\text{c}(s) =
       \frac12
       (s-1)^2
       +
       \frac{l_\text{c}^2}2|\nabla s|^2
       \,,
         \label{eq:surface}
     \end{align}
depending on the Griffiths constant $G_\text{c} >0$ and
a small length scale parameter $l_\text{c}>0$ which determines
the width of the diffusive interface \cite{Pandolfi_etal2021}.

\clearpage

\Paragraph{Phase-field evolution}
The initiation and propagation of phase-field fracture is modeled by the following principles:
\begin{itemize}
\item
The crack evolution depends on a retardation time $\tau_\text{r}>0$ and a
driving force $Y$, i.e.,
    \begin{align}
      \label{eq:pf_dgl}
      \tau_\text{r} \dot s &= -Y(\vec \sigma,s)\qquad\text{ in } (0,T) \times \Omega
      \,.
    \end{align}
\item
The driving force is composed of the elastic driving force $Y_\text{el}$
and a geometric term $Y_\text{geom}$ corresponding to the crack surface energy
    \begin{align*}
      Y(\vec \sigma,s) &=
      \begin{cases}
        Y_\text{el}(\vec\sigma) - M_\text{geom}Y_\text{geom}(s) &
        s > 0\,,
        \\
        0 & s = 0\,,
        \end{cases}
    \end{align*}
and complemented by the admissibility conditions which enforces
$s(t,\vec x)\geq 0$.  \gruen{This includes a scaling parameter
$M_\text{geom}>0$ which is required for the calibration of the two components
$Y_\text{el}$ and $Y_\text{geom}$ and which depends on the formulation of the crack driving force.
For a comparison of different models for~$Y_\text{el}$ we refer to  \cite{bilgen2019crack}.
 }

\item
Here, we use the stress-based criterion for the elastic driving force
    \begin{align}
      Y_\text{el}(\vec \sigma) &=
      \max\Big\{\frac{\sigma_\text{I}}{\sigma_\text{c}} -1,0\Big\}
    \end{align}
which is active if the maximum principal stress $\sigma_\text{I}$
exceeds the resistance value $\sigma_\text{c} > 0$.
\item
The geometric term $Y_\text{geom}$ is the gradient flow with respect to $\gamma_\text{c}$, i.e.,
    \begin{align}
      Y_\text{geom}(s)
      &=
      1-s
      +
      l_\text{c}^2\Delta s
      \,,
    \end{align}
complemented with homogeneous Neumann boundary conditions
$\nabla s\cdot\vec n = 0$ on $\partial\Omega$. The latter is realized in weak form, i.e.,
    \begin{align*}
      \int_\Omega Y_\text{geom}(s)\,\phi\,\mathrm d\vec x
      &=
      \int_\Omega
      \Big(
      (1-s)\phi
      -
      l_\text{c}^2\nabla s \cdot \nabla \phi\Big)\,\mathrm d\vec x
      \,,\qquad
      \phi\in \mathrm H^1(\Omega)
      \,.
    \end{align*}
\end{itemize}
Altogether, we obtain for the phase-field evolution \eqref{eq:pf_dgl} the parabolic equation
    \begin{align}
      \label{eq:pf_dgl_var}
      \tau_\text{r}
      \int_\Omega \dot s \phi\,\mathrm d\vec x
      +
      M_\text{geom}
      \int_\Omega
      \Big((s-1)\phi
      +
      l_\text{c}^2
      \nabla s \cdot \nabla \phi
      \Big)
      \,\mathrm d\vec x
      &=
      -
      \int_\Omega
      Y_\text{el}(\vec \sigma)\phi
      \,\mathrm d\vec x
      \,,\qquad
      \phi\in \mathrm H^1(\Omega)
      \,,
    \end{align}
which is complemented by a projection to ensure admissibiliy $s(\vec x) \in [0,1]$.

\Paragraph{Material degradation}

If $s(t,\vec x)\leq s_\text{min}$ at a material point $\vec x\in
\Omega$, then the material is cracked.  Depending on the phase-field
history $s_\text{inf}(t, \vec x)$ and the elastic domain
$\Omega\left(s_\text{inf}\right)$ the elasticity tensor at time $t$ is
defined by
    \begin{align}
      \label{eq:deg}
      \vec C(t,\vec x)
       &=
       s_\text{inf}(t',\vec x)
       \vec C
       +
       \big(1 - s_\text{inf}(t',\vec x)\big)
       \vec C_\text{reg}
       \,,
       \qquad
       t' = \inf \Big\{\tau\in [0,t]\colon
       \Omega\big(s_\text{inf}(\tau)\big) = \Omega\left(s_\text{inf}\right)\Big\}
    \end{align}
with a small but positive definite tensor $\vec C_\text{reg}$,
\gruen{e.g., $\vec C_\text{reg} = 10^{-7} \vec C$.}
\gruen{In the phase-field approximation, $\vec C_\text{reg}^{-1}$ is a penalty term
which enforces that the stress nearly vanishes in the crack zone.
Numerically we observe that the results are not very sensitive with respect to choice
of $\vec C_\text{reg}$.}

By this construction the material only degrades when the fracture
zone (determined by the parameter $s_\text{min}$) increases. In
particular, the full domain is elastic, i.e., $\vec C(0) = \vec C(t)$
as long $\Omega = \Omega\left(s_\text{inf}\right)$.

Then, the linear wave equation \eqref{eq:elastic} in $\Omega \setminus\Gamma_\text{c}(t)$
is approximated by
      \begin{align}
      \label{eq:elastic_degrad}
      \varrho \bdot {\vec v}(t) = \div   \vec \sigma(t)+ \vec f(t)
      \,,\qquad
      \bdot {\vec\sigma}(t) =
      \vec C(t)\vec \varepsilon(\vec v(t))
      \qquad\text{in }\Omega
      \,,
      \end{align}
and the phase-field approximation of the \gruen{mechanical energy} \eqref{eq:ee} is given by
    \begin{align*}
      \mathcal E_\text{pf}(t) =
      \frac12
       \int_\Omega
       \Big(\varrho|\vec v(t)|^2
       +\vec\sigma(t)\cdot\vec C(t)^{-1}\vec \sigma(t)\Big)\,\mathrm d\vec x
       \,.
    \end{align*}
This approximates the energy dissipation
     \begin{align*}
       \mathcal D_\text{pf}(t_1,t_2)
       =
       \mathcal E_\text{pf}(t_1) -  \mathcal E_\text{pf}(t_2)
       \geq 0
       \,,\qquad t_1<t_2\,,
     \end{align*}
which is strictly positive in case of fracture evolution, i.e.,  $\vec C(t_1) \neq \vec C(t_2)$.

Note that the degradation \eqref{eq:deg} is the most simple choice, in more general
formulations a monotone increasing degradation function
$g(\cdot)$ with $g(0) = 0$ and $g(1) = 1$ is included in \eqref{eq:deg};
see, e.g., \cite{kuhn2015degradation,sargado2018high}
for comparing different options.
Then, the degraded elasticity tensor is given by
    \begin{align*}
       \vec C_g(t,\vec x)
       &=
       g\big(s_\text{inf}(t',\vec x)\big)
       \vec C
       +
       \Big(1 - g\big(s_\text{inf}(t',\vec x)\big)\Big)
       \vec C_\text{reg}
       \,.
    \end{align*}
This can be complemented by a degradation of the density
    \begin{align*}
      \varrho_g(t,\vec x)
       &=
       g\big(s_\text{inf}(t',\vec x)\big)
       \varrho
       +
       \Big(1 - g\big(s_\text{inf}(t',\vec x)\big)\Big)
       \varrho_\text{reg}
    \end{align*}
and a moderate choice $\varrho_\text{reg} \in (0,\varrho]$.
The limiting physical model for $\ell_\text{c}\longrightarrow 0$ and
$\vec C_\text{reg}\longrightarrow 0$ with fixed $\varrho_\text{reg}>0$
is identical for all choices of material degradation.  For numerical
tests with sufficiently small length scale parameters $l_\text{c}>0$
and time steps ${\vartriangle} t>0$, see Section~\ref{sec:numerics},
no qualitative improvement is observed by including different degradation
functions, and so we formulate the algorithmic approach for the simple
choice~\eqref{eq:deg}.


\section{\gruen{A Runge-Kutta discontinuous Galerkin method for elastic waves}}\label{sec:DG}

The linear wave equation \eqref{eq:elastic_degrad} with time-dependent
elasticity tensor $\vec C(t)$ is approximated with a discontinuous
Galerkin (DG) method \gruen{in space and a Runge-Kutta method in
time}, see \cite{HocPSTW15} for details. \gruen{For hyperbolic
applications, this scheme is a widely used extension of finite volume
methods, see, e.g., \cite{dumbser2006arbitrary}. A first-order DG scheme is less
dissipative than a second-order approach and it provides also
approximations in case of discontinuities of the solution.}  Here we
shortly summarize the adoption of this method to our phase-field
fracture application.

The discretization is based on
\gruen{a formulation of the linear first-order system \eqref{eq:elastic_degrad}}
    \begin{align}
      \label{eq:system}
      M(t)
      \begin{pmatrix} \bdot {\vec v}(t) \\ \bdot {\vec\sigma}(t)\end{pmatrix}
      =
      A
      \begin{pmatrix}  {\vec v}(t) \\  {\vec\sigma}(t)\end{pmatrix}
       + \vec b(t)
        \,,\qquad t\in (0,T)
    \end{align}
\gruen{where the energy operator $M(t)$,
the differential operator $A$, and the right-hand side $\vec b$ are given by}
    \begin{align*}
      \big(M(t) (\bdot {\vec v}(t),\bdot {\vec\sigma}(t),(\vec w,\vec \eta)\big)_{0,\Omega}
      &=
      \big(\rho \bdot {\vec v}(t),\vec w(t)\big)_{0,\Omega}
      +
      \big(\vec C(t)^{-1}\bdot {\vec \sigma}(t),\vec \eta\big)_{0,\Omega}
      \,,
      \\
      \big(A(\vec v(t),\vec\sigma(t)),(\vec w,\vec \eta)\big)_{0,\Omega}
      &=
      \big(\div {\vec \sigma}(t),\vec w\big)_{0,\Omega}
      +
      \big(\vec\varepsilon({\vec v}(t),\vec \eta\big)_{0,\Omega}
     +
     \big(\vec v(t),\vec \eta\vec n\big)_{0,\Gamma_\text{D}}
     +
     \big( \vec\sigma(t)\vec n,\vec w\big)_{0,\Gamma_\text{N}}
      \\
     \big(\vec b(t),(\vec w,\vec \eta)\big)_{0,\Omega}
     &
     =
     \big(\vec f(t),\vec w\big)_{0,\Omega}
     +
     \big(\vec v_\text{D}(t),\vec \eta\vec n \big)_{0,\Gamma_\text{D}}
     +
     \big(\vec g_\text{N}(t),\vec w\big)_{0,\Gamma_\text{N}}
    \end{align*}
\gruen{for test functions $(\vec w,\vec \eta)$ in $\Omega$. This defines $
      \mathcal E_\text{pf}(t) =
      \frac12 \big(M(t) ({\vec v}(t),{\vec\sigma}(t),({\vec v}(t),{\vec\sigma}(t))\big)_{0,\Omega}
      $.
We use the standard notation for
the $\mathrm L_2$ inner products
$(\cdot,\cdot)_{0,\Omega}$ and $(\cdot,\cdot)_{0,\Gamma}$
in the domain and on the boundaries.}

For the approximation of \eqref{eq:system} in space, we need to construct approximations
$M_h$, $A_h$ and $\vec b_h$.
On a mesh $\Omega_h = \bigcup_{K\in \mathcal K_h} K$ with
elements $K$, let $V_h^\text{dg} = \prod_{K\in \mathcal K_h} \mathbb
P_k(K;\Real^d\times \Real^{d\times d}_\text{sym})$ be the
discontinuous finite element space of polynomial degree $k$.
\gruen{The discrete energy operator
$M_h(t) = \sum_{K\in \mathcal K_h} M_{h,K}(t)$
for the discontinuous functions $\vec v_h = \sum_{K\in \mathcal K_h} \vec v_{h,K}$ and
$\vec \sigma_h = \sum_{K\in \mathcal K_h} \vec \sigma_{h,K}$ is defined locally on $K$ by}
    \begin{align*}
      \big(M_{h,K}(t) ({\vec v}_{h,K}, {\vec\sigma}_{h,K}),(\vec w_{h,K},\vec \eta_{h,K})\big)_{0,K}
      &=
      \big(\rho {\vec v}_{h,K},\vec w_h\big)_{0,K}
      +
      \big(\vec C(t)^{-1} {\vec \sigma}_{h,K},\vec \eta_h\big)_{0,K}
      \,.
    \end{align*}
For the discontinuous functions, the derivatives are approximated by jump
terms on the faces $\mathcal F = \bigcup_K \mathcal F_K$,
where $\mathcal F_K$ are the faces on every element $K$.  For inner
faces $f \in \mathcal F\cap \Omega$, let $K_f$ be the neighboring cell
such that $\bar f= \partial K \cap \partial K_f$. On boundary faces $f
\in \mathcal F\cap \partial\Omega$ we set $K_f = K$.  Let $\vec n_K$
be the outer unit normal vector on~$\partial K$.  We define the jump
$[\vec v_h]_{K,f} = \vec v_{h,K_f} -\vec v_{h,K}$ on inner faces,
where $\vec v_{h,K}$ denotes the continuous extension of $\vec v_h|_K$
to $\overline K$. In the same way, the jump for the stress tensor is defined.
On Dirichlet boundary faces, we set $[\vec v_h]_{K,f} = \vec 0$
and $ [\vec \sigma_h]_{K,f}\vec n = -2\, \vec \sigma_h\vec n$.
On Neumann boundaries, set $[\vec v_h]_{K,f} =-2\vec v_h$
and $[\vec \sigma_h]_{K,f}\vec n = \vec 0$.

The full upwind DG appro\-ximation $A_h(t) = \sum_{K\in \mathcal K_h} A_{h,K}(t)$
is defined by the local contributions 
    \begin{align*}
      &\big(A_{h,K}(t) (\vec v_h,\vec\sigma_h),(\vec w_h,\vec \eta_h)\big)_{0,K}
      =
      \big(\div \vec \sigma_{h,K},\vec \psi_{h,K}\big)_{0,K}
      +
      \big(\vec \varepsilon(\vec v_{h,K}), \vec\eta_{h,K}\big)_{0,K}
      \\
      &\qquad
      +
      \frac12
      \sum_{f\in \mathcal F_K}
      \Big(\vec n_K\cdot \big([\vec\sigma_{h,K}]_{K,f}\vec n_K
        +Z_{\text{P}}(t)[\vec v_h]_{K,f}\big)
        ,
        \vec n_K\cdot \big(Z_{\text{P}}(t)^{-1} \vec \eta_{h,K}\vec n_K  +
         \vec w_{h,K}\big)\Big)_{0,f}
      \\
      &\qquad
      +
      \frac12
      \sum_{f\in \mathcal F_K}
      \Big(\vec n_K\times \big([\vec\sigma_{h,K}]_{K,f}\vec n_K
        +Z_{\text{S}}(t)[\vec v_h]_{K,f}\big)
        ,
        \vec n_K\times \big(Z_{\text{S}}(t)^{-1} \vec \eta_{h,K}\vec n_K  +
        \vec w_{h,K}\big)\Big)_{0,f}
    \end{align*}
depending on the impedances $Z_{\text{P}}(t) =
\sqrt{\varrho(2\mu(t)+\lambda(t))}$ and $Z_{\text{S}}(t) =
\sqrt{\varrho\mu(t)}$ of compressional waves and shear waves,
respectively.  Here, $(\cdot,\cdot)_{0,K}$ and $(\cdot,\cdot)_{0,f}$
denotes the $\mathrm L_2$ inner product in the elements and on the
faces, respectively. The material parameters and thus the upwind flux
and the operator $A_h$ depend on the material degradation
\eqref{eq:deg} encoded in $\vec C(t)$.
Note that the construction of the upwind flux can be extended
to composite materials and discontinuous material
parameters since it is computed by the exact solution of Riemann problems at interfaces.
Here, the material tensor $\vec C(t)$ is continuous in space which
simplifies the evaluation of the numerical flux on the element faces.

The boundary and volume data enter in the right-hand side
$\vec b_h(t) = \sum\nolimits_{K\in \mathcal K_h} \vec b_{h,K}(t)$ with
    \begin{align*}
     \big(\vec b_{h,K}(t),(\vec w_{h,K},\vec \eta_{h,K})\big)_{0,K}
     \!
     =
     \big(\vec f(t),\vec w_{h,K}\big)_{0,K}
     &
     +
     \!
     \!
     \!
     \!
     \sum_{f\in \mathcal F_K\cap \partial_\text{D}\Omega}
     \Big(
      \big(\vec n_K\cdot \vec v_\text{D}(t)
        ,
        \vec n_K\cdot (Z_{\text{P}}(t)^{-1} \vec \eta_{h,K}\vec n_K  +
        \vec w_{h,K})\big)_{0,f}
        \\[-1mm]
        &\qquad\qquad\quad
        +
      \big(\vec n_K\times \vec v_\text{D}(t)
        ,
        \vec n_K\times (Z_{\text{S}}(t)^{-1} \vec \eta_{h,K}\vec n_K  +
        \vec w_{h,K})\big)_{0,f}
        \Big)
      \\
      &
      +
     \!
     \!
     \!
     \!
      \sum_{f\in \mathcal F_K\cap \partial_\text{N}\Omega}
      \Big(
      \big(\vec n_K\cdot \vec g_\text{N}(t)
        ,
        \vec n_K\cdot (\vec \eta_{h,K}\vec n_K  +
        Z_{\text{P}}(t) \vec w_{h,K})\big)_{0,f}
        \\[-1mm]
        &\qquad\qquad\quad
        +
      \big(\vec n_K\times \vec g_\text{N}(t)
        ,
        \vec n_K\times (\vec \eta_{h,K}\vec n_K  +
        Z_{\text{S}}(t) \vec w_{h,K})\big)_{0,f}
        \Big)
      \,.
   \end{align*}
For the discretization in time, we distinguish two cases.
\begin{description}
  \item[Elastic time step]
If $\vec C(t_n) = \vec C(t_{n-1})$ in the time step from $t_{n-1}$ to $t_n$,
we use the implicit midpoint rule with time step size ${\vartriangle} t_n = t_n - t_{n-1}$. This corresponds for 
$(\vec v_h^n,\vec \sigma_h^n)$ to the equation
    \begin{align*}
      \frac1{{\vartriangle} t_n} M_h^{n-1}
      \gruen{\begin{pmatrix}  {\vec v}^n -{\vec v}^{n-1} \\  {\vec\sigma}^n-{\vec\sigma}^{n-1}\end{pmatrix}}
      = \frac12 A_h^{n-1}
      \begin{pmatrix}  {\vec v}^n +{\vec v}^{n-1} \\  {\vec\sigma}^n+{\vec\sigma}^{n-1}\end{pmatrix}
      + \vec b_h^{n-1/2}
      \,,
    \end{align*}
with $M_h^{n-1}= M_h(t_{n-1})$, $A_h^{n-1}= A_h(t_{n-1})$, and $\vec b^{n-1/2}_h =
\vec b_h\big(\frac12(t_n+t_{n-1})\big)$,
so that 
$(\vec v_h^n,\vec \sigma_h^n)$ is determined by solving the linear system
    \begin{align}
      \label{eq:midpoint}
      \Big( M_h^{n-1} - \frac{{\vartriangle} t_n}2 A_h^{n-1} \Big)
      \begin{pmatrix}  {\vec v}^n  \\  {\vec\sigma}^n\end{pmatrix}
          =
          \Big( M_h^{n-1} + \frac{{\vartriangle} t_n}2 A_h^{n-1}\Big)
      \begin{pmatrix}  {\vec v}^{n-1} \\  {\vec\sigma}^{n-1}\end{pmatrix}
          +
              {\vartriangle} t_n\,
              \vec b_h^{n-1/2}
              \,.
    \end{align}
This scheme is energy conserving for conforming solutions,
so that in case of $\vec b_h^{n-1/2} = \vec 0$
the \gruen{mechanical energy} is conserved
$\mathcal E_\text{pf}(t_n) \approx  \mathcal E_\text{pf}(t_{n-1})$
up to the numerical dissipation induced by the nonconforming DG approximation.
  \item[Dissipative time step]
If the material degrades in the time step, i.e., $\vec C(t_n) \neq \vec C(t_{n-1})$, we use the implicit Euler method
    \begin{align*}
      M_h^n
            \begin{pmatrix}  {\vec v}^n  \\  {\vec\sigma}^n\end{pmatrix} =
              M_h^{n-1}
                    \begin{pmatrix}  {\vec v}^{n-1} \\  {\vec\sigma}^{n-1}\end{pmatrix}
      +
      {\vartriangle} t_n
      A_h^n
            \begin{pmatrix}  {\vec v}^n  \\  {\vec\sigma}^n\end{pmatrix}
          +
              {\vartriangle} t_n\,
              \vec b_h^n
    \end{align*}
with $M_h^n= M_h(t_n)$, $A_h^n= A_h(t_n)$, and $\vec b^n_h =
\vec b_h(t_n)$,
so that $\vec y_h^n$ is determined by the linear system
    \begin{align}
      \label{eq:Euler}
      \Big( M_h^n - {\vartriangle} t_n A_h^n \Big)
      \begin{pmatrix}  {\vec v}^n  \\  {\vec\sigma}^n\end{pmatrix}
          =
           M_h^{n-1}
           \begin{pmatrix}  {\vec v}^{n-1} \\  {\vec\sigma}^{n-1}\end{pmatrix}
           +
              {\vartriangle} t_n\,
              \vec b_h^n
              \,.
    \end{align}
For $\vec b_h^{n} = \vec 0$ the difference $\mathcal E_\text{pf}(t_{n-1}) -\mathcal E_\text{pf}(t_n)$
corresponds to the dissipated energy by opening the fracture zone so that
the elastic domain gets smaller, i.e.,
$\Omega\big(s_\text{inf}(t_n)\big)\subset\Omega\big(s_\text{inf}(t_{n-1})\big)$.
\end{description}


The linear systems \eqref{eq:midpoint} and \eqref{eq:Euler} are
well-defined, since the upwind discretization $-A_h(t)$ is positive
semi-definite, so that the matrices $M_h^{n-1} - \frac{{\vartriangle}
  t_n}2 A_h^{n-1}$ and $M_h^n - {\vartriangle} t_n A_h^n $ are
regular.  Moreover, both schemes are unconditional stable, also if the
\gruen{Courant-Friedrichs-Lewy} 
condition $c_\text{max}{\vartriangle} t_n\leq \frac 12 h$ is not
satisfied; herein, the maximal wave speed is denoted by
$c_\text{max} = \max\sqrt{(2\mu+\lambda)/\varrho}$.  On the other hand,
$c_\text{max}{\vartriangle} t_n = \mathcal O(h)$ is required for a
balanced approximation error in space and time. Then we observe that
the linear systems in \eqref{eq:midpoint} and \eqref{eq:Euler} are
well conditioned \cite[Lem.~3.2]{Rieder2020}, so that the approximate
solution can be computed iteratively
within a small number of steps, e.g., using a parallel GMRES
method with block-Jacobi preconditioning.

In the fracture zone where the material is degraded, the material
stiffness and thus also the wave speed is reduced. In the limit this
approximates Neumann boundary conditions at crack interfaces.

\clearpage


\section{Computational dynamic fracture}\label{sec:fracure}

We combine the DG discretization for the hyperbolic linear wave
equation in the discontinuous finite element space
$V_h^\text{dg}\subset \mathrm L_2(\Omega;\Real^d\times \Real^{d\times
  d}_\text{sym})$ with a parabolic phase-field approximation by lowest
order conforming finite elements $V_h^\text{cf}\subset \prod \mathbb
P(K)\cap\mathrm C(\bar\Omega)\subset \mathrm H^1(\Omega)$.  Then,
$\phi_h\in V_h^\text{cf}$ is uniquely defined by the nodal values
$\big(\phi_h(\vec x)\big)_{\vec x\in \mathcal C_h}$, where $\mathcal
C_h\subset \overline\Omega$ 
are the vertices of the elements $K\in\mathcal K_h$.
Since the finite element approximation of the phase field is
continuous, it allows a continuous evaluation of the degraded material parameters on the element faces within the computation of the
numerical flux in the wave discretization.

\begin{figure}[H] 
  \centering
\begin{tikzpicture}[node distance=2.5cm]
  \node (inp)
	[startstop]
	{\parbox{7.75cm}{\centering
	    set initial values $(\vec v_h^0,\vec \sigma_h^0)$\\
            start with elastic domain $\Omega^0= \Omega$ and $s_h^0=1$}
  };
  \node (pro1)
	[process, below of=inp, yshift=0.5cm]
	{\parbox{10.5cm}{\centering
            compute $(\vec v^n_h,\vec \sigma_h^n)$ from $(\vec v_h^{n-1},\vec \sigma_h^{n-1})$
            with implicit midpoint rule
            \\
             using the material parameters from the previous time step
        }};
  \node (pro2)
	[process, below of=pro1, yshift=1cm]
	{\parbox{7cm}{\centering
            compute $s_h^n$ from $\vec \sigma_h^n$ and $s_h^{n-1}$
        }};
  \node (pro3) [process, below of=pro2, yshift=0.8cm, fill=lightgreen]
	{\parbox{7cm}{\centering
            project $s_h^n(\vec x)$ to the admissible range $[0,1]$\\
            set $s_h^n(\vec x)$ to zero for $s_h^n(\vec x) < s_\text{min}$
            \\
            compute the elastic domain $\Omega^n$
	}};
  \node (pro4) [process, left of=pro3, xshift=-1.0cm, yshift=-2.5cm, fill=lightgreen]
		{\parbox{3.5cm}{\centering
		    the time step is elastic
	}};
  \node (pro5) [process, right of=pro3,xshift=0.0cm,yshift=-2.5cm, fill=lightgreen]
		{\parbox{4.75cm}{\centering
                    the fracture zone is growing
	}};
  \node (pro51) [process, below of=pro5,xshift=0.0cm,yshift=1.2cm, fill=lightgreen]
		{\parbox{9.5cm}{\centering
                    update the material parameters depending on $s_{h,\text{inf}}^n$
	}};
  \node (pro6) [process, below of=pro51, yshift=1cm, xshift=0.0cm]
		{\parbox{9.5cm}{\centering
                    repeat the computation of
                    $(\vec v^n_h,\vec \sigma_h^n)$,
                    now with the updated material parameter and using the implicit Euler method	
	}};
	\draw [arrow] (inp.south) -- node[left] {$n:=1$} (pro1.north);
	\draw [arrow] (pro1.south) -- node[left] {} (pro2.north);
	\draw [arrow] (pro2.south) -- node[left] {} (pro3.north);
	\draw [arrow] (pro3.south) -- +(0,-0.75) -| node[above left]{$\Omega^{n-1}\setminus \Omega^n\neq \emptyset$} (pro5.north) ;
	\draw [arrow] (pro3.south) -- +(0,-0.75) -| node[above right]{$\Omega^n= \Omega^{n-1}$} (pro4.north);
	\draw [arrow] (pro4.west) -- +(-2.5,0) |- node[above right]{$n := n + 1$} (pro1.west);
	\draw [arrow] (pro5.south) -- (pro51.north);
	\draw [arrow] (pro51.south) -- (pro6.north);
	
	\draw [arrow] (pro6.east) -- +(0.5,0) |- node[above left]{$n := n + 1$} (pro1.east);
\end{tikzpicture}
\caption{The staggered scheme for
  the elastic variables $(\vec v_h^n,\vec \sigma_h^n)$ and the phase field $s_h^n$.}
\label{fig:staggered}
\end{figure}
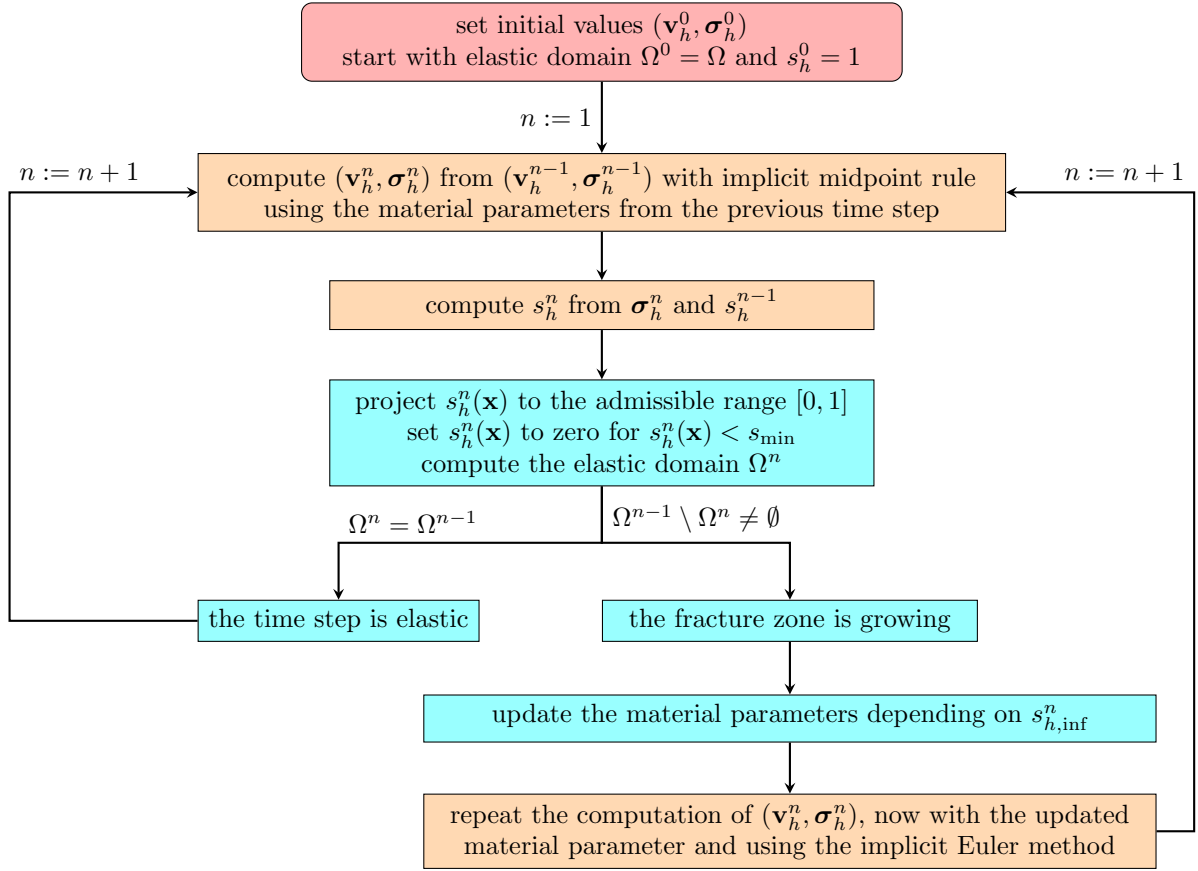

We use a staggered time-discrete scheme, updating alternately the
elastic system $(\vec v_h^n,\vec \sigma_h^n)$ and the phase field~$s_h^n$,
and depending on the phase-field evolution the material is
degraded only in time steps where the elastic domain (determined by
$s_\text{min}$) becomes smaller. The main steps of the algorithm are
summarized in Fig.~\ref{fig:staggered}.

\bigskip

\clearpage

We start with initial values $(\vec v_h^0,\vec \sigma_h^0,s_h^0) \in
V_h^\text{dg}\times V_h^\text{cf}$ in the material without fracture,
i.e., we set $\Omega^0= \Omega$, $s_h^0 = s_{h,\text{inf}}^0= 1$, and
$\vec C^0 = \vec C$.  For elasticity we set $\vec v_h^0= \Pi_h\vec
v_0$ and $\vec \sigma_h^0= \Pi_h\vec \sigma_0$ with a suitable
projection or interpolation $\Pi_h$.  Furthermore, we select time
steps $0<{\vartriangle} t_\text{pf}\leq {\vartriangle} t_\text{el}$
for the evolution of the phase-field and the elastic wave propagation.
In addition, the displacement can be approximated by integrating the
velocity; we set $\vec u_h^0= \Pi_h\vec u_0$.
We start with ${\vartriangle} t_1 ={\vartriangle} t_\text{el}$.

In every time step $n=1,2,3,\ldots$ we proceed as follows:
\begin{itemize}
\item[(S1)]
We compute from
$(\vec v^{n-1}_h,\vec \sigma^{n-1}_h)$ the solution for the
next time step $(\vec v^n_h,\vec \sigma^n_h)\in V_h^\text{dg}$
with the implicit midpoint rule \eqref{eq:midpoint}, i.e.,
    \begin{align*}
      &\Big(\big(M_h^{n-1}-\frac{{\vartriangle} t_n}2 A_h^{n-1}\big)
      (\vec v_h^n,\vec\sigma_h^n),(\vec w_h,\vec \eta_h)\Big)_{0,\Omega}
      \\
      &\quad \quad
      =
      \Big(\big(M_h^{n-1}+\frac{{\vartriangle} t_n}2 A_h^{n-1}\big)(\vec v_h^{n-1},\vec\sigma_h^{n-1})
      ,(\vec w_h,\vec \eta_h)\Big)_{0,\Omega}
      +
      {\vartriangle} t_n
      \big(\vec b_{h}^{n-1/2},(\vec w_h,\vec \eta_h)\big)_{0,\Omega}
      \,,\quad
      (\vec w_h,\vec \eta_h) \in V_h^\text{dg}
      \,.
    \end{align*}
Here, $(\vec v^n_h,\vec \sigma^n_h)$ is just a candidate for the next time step
which will be accepted
only if the fracture criterion in step (S3) shows that this step is elastic; otherwise,
$(\vec v^n_h,\vec \sigma^n_h)$ will be recomputed in (S6).
\item[(S2)] Depending on $\vec \sigma_h^n$, we approximate the phase field $s_h^n\in V_h^\text{cf}$
by the implicit Euler method applied to~\eqref{eq:pf_dgl_var}, i.e., by solving
    \begin{align*}
      \big(
      \tau_\text{r} s_h^n
      ,
      \phi_h
      \big)_{0,\Omega}
      +
      {\vartriangle} t_n
      M_\text{geom}
      \Big(
      \big(
      (s_h^n-1)
      ,
      \phi_h
      \big)_{0,\Omega}
      +
      \big(
      l_\text{c}^2
      \nabla s_h^n , \nabla \phi_h
      \big)_{0,\Omega}
      \Big)
      &=
      \big(
      \tau_\text{r} s_h^{n-1}
      -
      {\vartriangle} t_n
      Y_\text{el}(\vec \sigma_h^n)
      ,
      \phi_h
      \big)_{0,\Omega}
      \,,\quad \phi_h\in V_h^\text{cf}
      \,.
    \end{align*}
Again we note that $s_h^n$ is just a candidate since it will be modified in the next step (S3).
\item[(S3)]
On all nodal points $\vec x\in \mathcal C_h$, the phase field $s_h^n(\vec x)$ is projected to $[0,1]$ by
    \begin{align*}
     s_h^n(\vec x) &\coloneqq
     \begin{cases}
       1                             &
       s_h^n(\vec x)\geq 1\,,\\
       s_h^n(\vec x)
       &  s_\text{min} \leq s_h^n(\vec x) < 1\,,\\
       0            &   s_h^n(\vec x)< s_\text{min}
                        \text{ or } s^{n-1}(\vec x) = 0\,,
      \end{cases}
    \end{align*}
and we set
    \begin{align*}
      s_{h,\text{inf}}^n(\vec x) = \min\big\{s_{h,\text{inf}}^{n-1}(\vec x),s_h^n(\vec x)\big\}\,,
      \qquad
      \Omega^n = \Big\{\vec x\in \Omega^{n-1}\colon s_{h,\text{inf}}^n(\vec x)
      \geq s_\text{min}\Big\}\,.
    \end{align*}
\item[(S4)]
  If $\Omega^n= \Omega^{n-1}$, the time step is elastic and $\vec C^n = \vec C^{n-1}$;
  then, we accept the elastic solution
  $(\vec v_h^n,\vec \sigma_h^n)$ and   directly proceed with (S7).
\item[(S5)]
Otherwise, if the fracture zone is growing by $\Omega^{n-1}\setminus\Omega^n\neq\emptyset$,
we update the material, i.e.,
    \begin{align*}
      \vec C^n(\vec x)
       =
       s_\text{inf}^n(\vec x)
       \vec C
       +
       \big(1 - s_\text{inf}^n(\vec x)\big)
       \vec C_\text{reg}
      \,.
    \end{align*}
\item[(S6)]
We repeat the computation of $(\vec v^n,\vec \sigma^n)$ with the updated material:\\
Using the implicit Euler method~\eqref{eq:Euler},
we compute $(\vec v_h^n,\vec \sigma_h^n)\in V_h^\text{dg}$ by solving
    \begin{align*}
      &\Big(\big(M_h^n- {\vartriangle} t_n A_h^n\big)
      (\vec v_h^n,\vec\sigma_h^n),(\vec w_h,\vec \eta_h)\Big)_{0,\Omega}
      \\
      &\qquad \qquad
      =
      \big(M_h^{n-1}(\vec v_h^{n-1},\vec\sigma_h^{n-1})
      ,(\vec w_h,\vec \eta_h)\Big)_{0,\Omega}
      +
      {\vartriangle} t_n
      \big(\vec b_{h}^n,(\vec w_h,\vec \eta_h)\big)_{0,\Omega}
      \,,\qquad
      (\vec w_h,\vec \eta_h) \in V_h^\text{dg}
      \,.
    \end{align*}
\item[(S7)]
We set $\vec u_h^n = \vec u_h^{n-1} + {\vartriangle} t_n\, \vec v_h^n$.
\\
If $s_h^n = s_h^{n-1}$, we expect that the next time step will also be elastic and we set
${\vartriangle} t_{n+1}={\vartriangle} t_\text{el}$;
otherwise, we set ${\vartriangle} t_{n+1}={\vartriangle} t_\text{pf}$.
Then, we continue with the next time step $n:= n+1$ and we proceed with (S1).
\end{itemize}

\clearpage

\Clearpage

\section{Numerical experiments}\label{sec:numerics}

The numerical versatility of our staggered algorithm is illustrated
with examples in one, two and three dimensions.

\subsection{\gruen{Simulation of a Hopkinson bar experiment}}

\begin{figure}[H]
\includegraphics[width=0.99\textwidth]{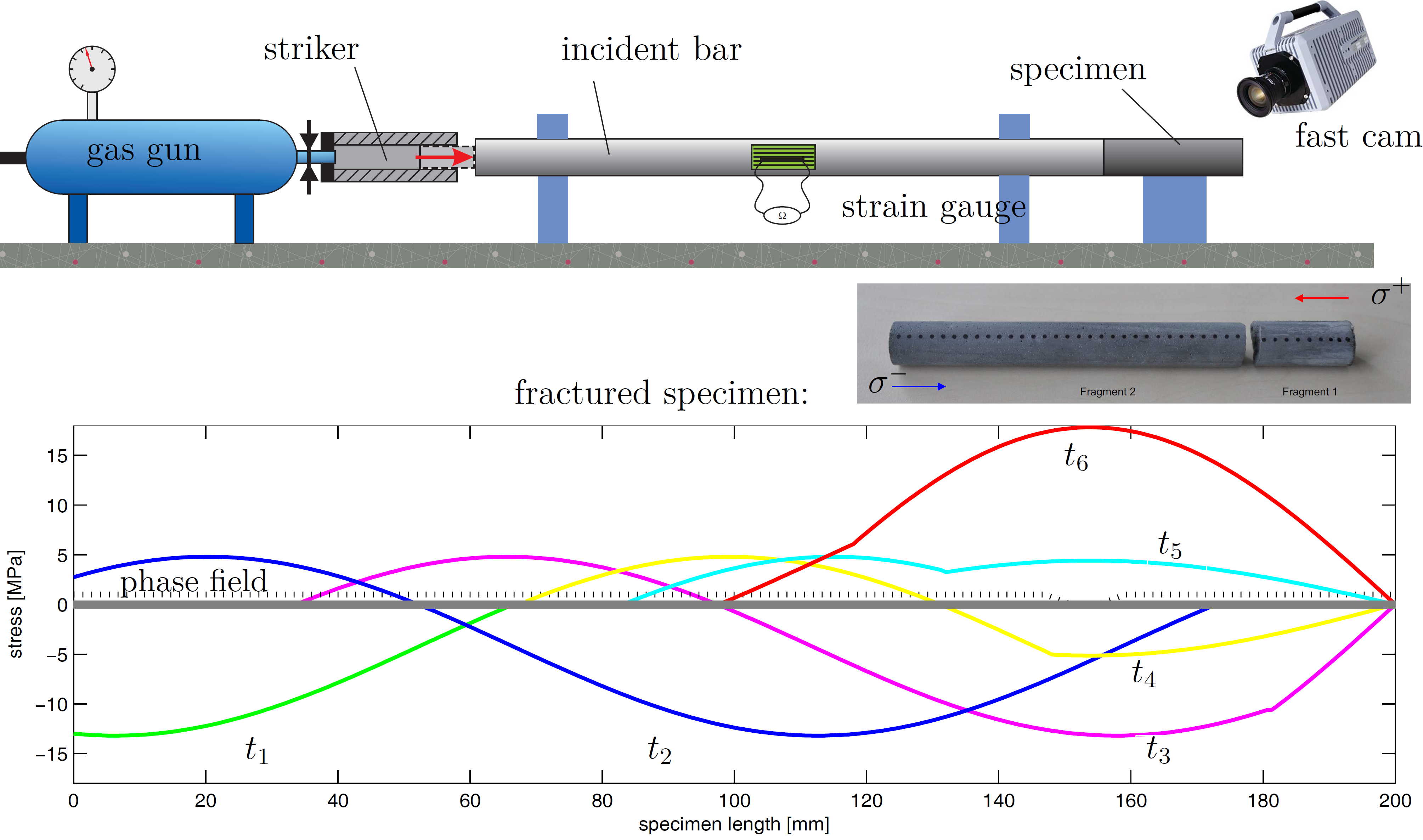}
 \caption{Hopkinson bar spallation experiment with a fractured specimen and simulation of the stress wave for different times; $t_1$: incoming pressure wave pulse, $t_2$: traveling wave, $t_3-t_4$: superposition during reflection at the free end, and $t_5$: the peak stress which causes spallation; the dashed line is the phase field. }   \label{fig:SHB}
\end{figure}

We start with the simulation of a spalling experiment performed in our
lab. In these experiments, the setup of a Hopkinson-Pressure bar is
modified so that the specimen is placed at the end of a long incident
bar that a striker hits.  The induced pressure pulse is transmitted
via the incident bar into the specimen and reflected at its free
end. The resulting tensional wave determines an inhomogeneous stress
state, and when the tensile stress exceeds the material's strength,
the specimen fractures, see Fig.~\ref{fig:SHB}.  The stress value at
fracture is considered to be the dynamic tensile resistance of the material,
cf. \cite{WeinbergKhosravani2018Dymat}. A successful experiment
requires specimens with a certain tension-compression asymmetry;
we investigated Ultra-High Performance Concrete  with $E=50\,$GPa
and $\sigma_c=18\,$MPa, cf. \cite{Khosravani_etmany2019UHPC}.

The evaluation of the experiment is based on the one-dimensional
linear wave theory and our DG simulation of this setting completely
reproduces the theoretical results. For some instances of time, the
waves in the specimen are displayed in Fig.~\ref{fig:SHB} whereby the
analytical and the numerical curve basically lay on top of each
other. The phase-field fracture approximation shows a crack at the
critical stress of 18\,MPa.

Clearly, this setting is too simple to illustrate the advantages of
our new numerical approach. Therefore we modified
it in such a way that the examples show waves propagation and
superposition in two or three dimensions. To avoid effects which are
only caused by symmetry, we choose a non-uniformly curved geometry
with non-symmetric wave pulses from left and right. The pulses travel
with the sample's wave speed, superpose, are then reflected at the
free boundary, and continue traveling with inverted amplitude until
they induce cracks. In the following, geometry and all material parameters are
dimensionless  but the values of choice (if understand in units N, mm
and \textmu s) correspond to a typical hard plastic such as polymethyl
methacrylate (PMMA).

\clearpage

\subsection{A 2D curved bar}

For plane strain computations a unit reference domain is mapped into a curved configuration $\Omega = \vec \varphi(\Omega_\text{ref})$
with
    \begin{align*}
      \Omega_\text{ref} = (-0.5,0.5)\times (0.03125,0.03125) 
      \ \text{ and }\
      \vec\varphi(x_1,x_2) = \big(x_1,\cos(0.5 x_1\pi)\big)
      +x_2\big(\sin(0.5 x_1\pi),\cos(0.5 x_1\pi)\big)\,.
    \end{align*}
\begin{figure}[H]
  \centering \vspace*{-4mm}
    \includegraphics[width=0.95\textwidth]{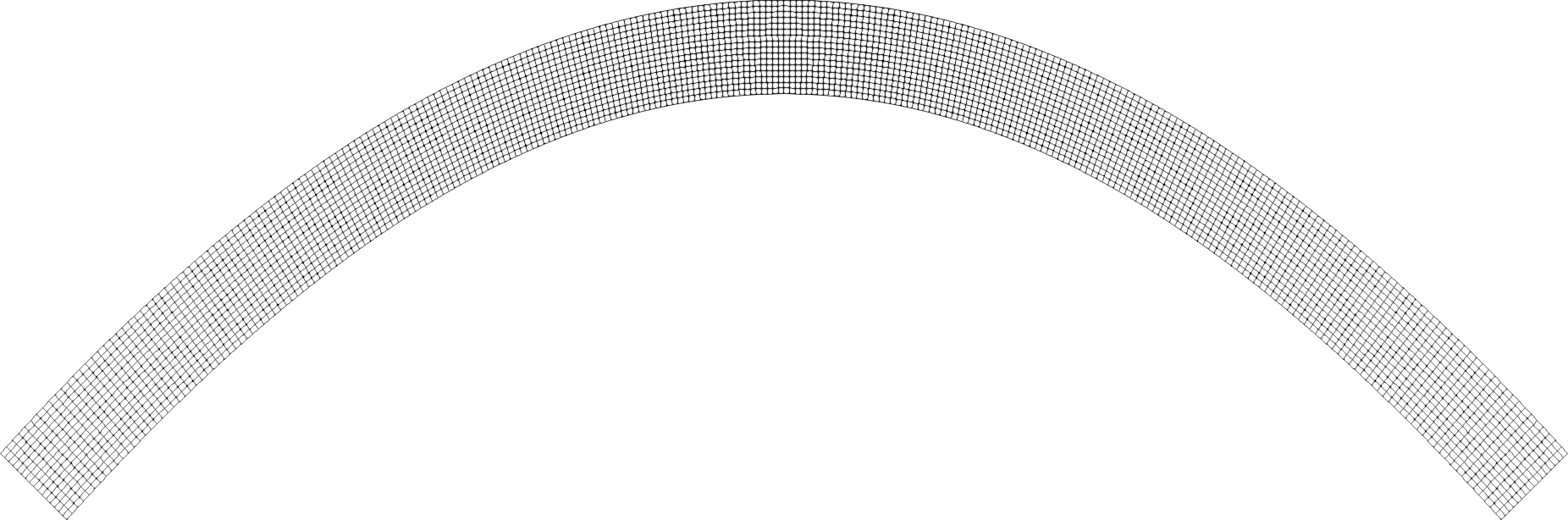} \\[-10mm]
    \caption{Coarse mesh  $\Omega_h\subset \Omega = \vec \varphi(\Omega_\text{ref})$\\
      with mesh size $h = 2^{-8}$ and 4096 quadrilaterals.}
\label{fig:mesh}
\end{figure}

The material is  isotropic  and linear elastic with
$\vec C\vec\varepsilon= 2\mu\vec\varepsilon+\lambda\trace(\vec\varepsilon)\id$,
$\mu = 1$,
$\lambda = 2$. With density $\rho = 1$ this gives the speed $c_\text{P}= 2$ for compressional waves
and $c_\text{S}= 1$ for shear waves.

For the phase-field model, we use
$\sigma_\text{c} = 27$,
$M_\text{geom} = 0.01$,
$l_\text{c} = 0.0005$,
$s_\text{min} = 0.01$,
and for the material degradation
$\vec C_\text{reg}\vec\varepsilon=
2\mu_\text{reg}\vec\varepsilon+\lambda_\text{reg}\trace(\vec\varepsilon)\id$
with $\mu_\text{reg} = 10^{-7}\mu$ and
$\lambda_\text{reg} = 10^{-7}\lambda$.
The phase-field is approximated with bilinear conforming finite
elements; the mesh is illustrated in~Fig.~\ref{fig:mesh} for the
coarsest level with mesh size $h = 2^{-8}$  and $16 \cdot 256 $ quadrilaterals, and
uniformly refined up to $h = 2^{-11}$ with $128 \cdot 2048 $ quadrilaterals.
Velocity and stress are approximated with discontinuous bilinear finite
elements.

The model is loaded by a smooth pressure pulse $\vec \sigma \vec n = g_\text{N}(t)\vec n $
at the left and right boundary with
    \begin{align*}
      g_\text{N}(t,\vec\varphi(\pm 0.5,x_2)) = a_\pm(c_\text{P}t-S_\pm)\,,
      \qquad
      a_\pm(s) = A_\pm\exp\Big(\frac{-1}{w_\pm^2-s^2}\Big) \,,
      \quad t\in (0,t_\text{init})
      \,.
    \end{align*}
The impulse width is $w_\pm = 0.3$ and the duration is $t_\text{init} = 0.24$.
This corresponds to incoming  compressional waves from left and right. To break symmetry the impulse on the right side is 5\% stronger, i.e.,
we set for the amplitude $A_+= 1.05 A_- $ 
and use the time shift parameters $S_+ = 1.25$ and $S_- = -1.03$.
In that way 
a crack will be initiated by the superposition of the traveling waves close to the center at $x_1= 0$ but not exactly aligned with the mesh.
On the remaining boundaries and for $t\geq t_\text{init}$ we use homogeneous Neumann boundary conditions $\vec g_\text{N}=0$.

\begin{figure}[H]
\begin{tabular}{c}
\includegraphics[width=0.47\textwidth]{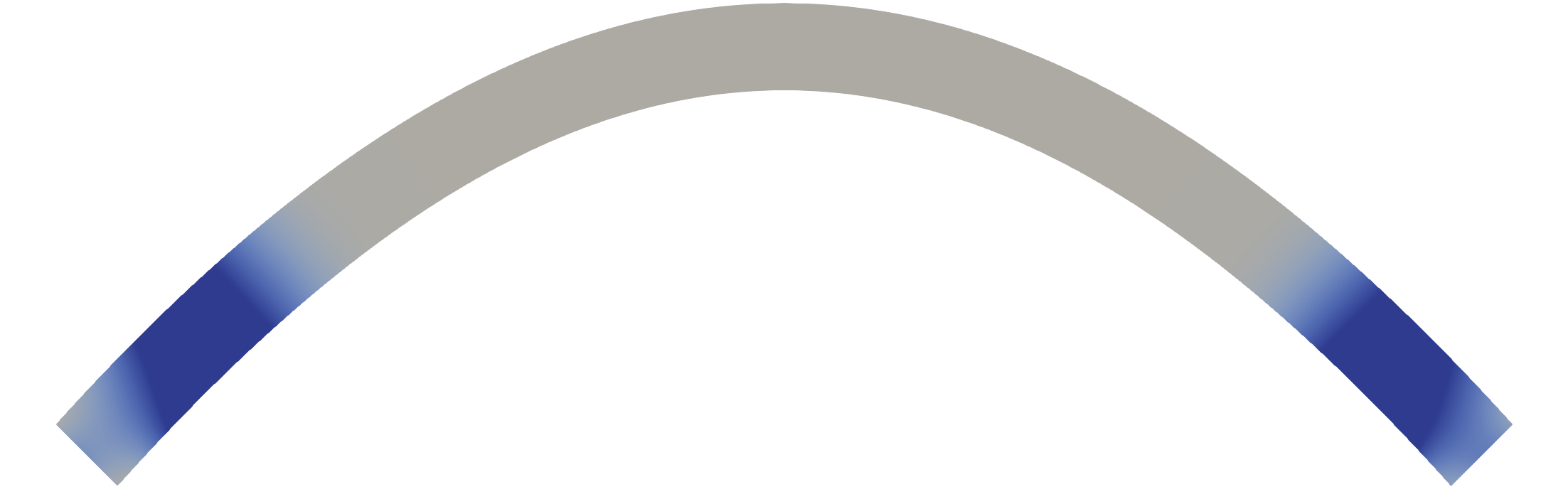}
\\[-1mm]
$t = 0.18$
\\[5mm]
\includegraphics[width=0.47\textwidth]{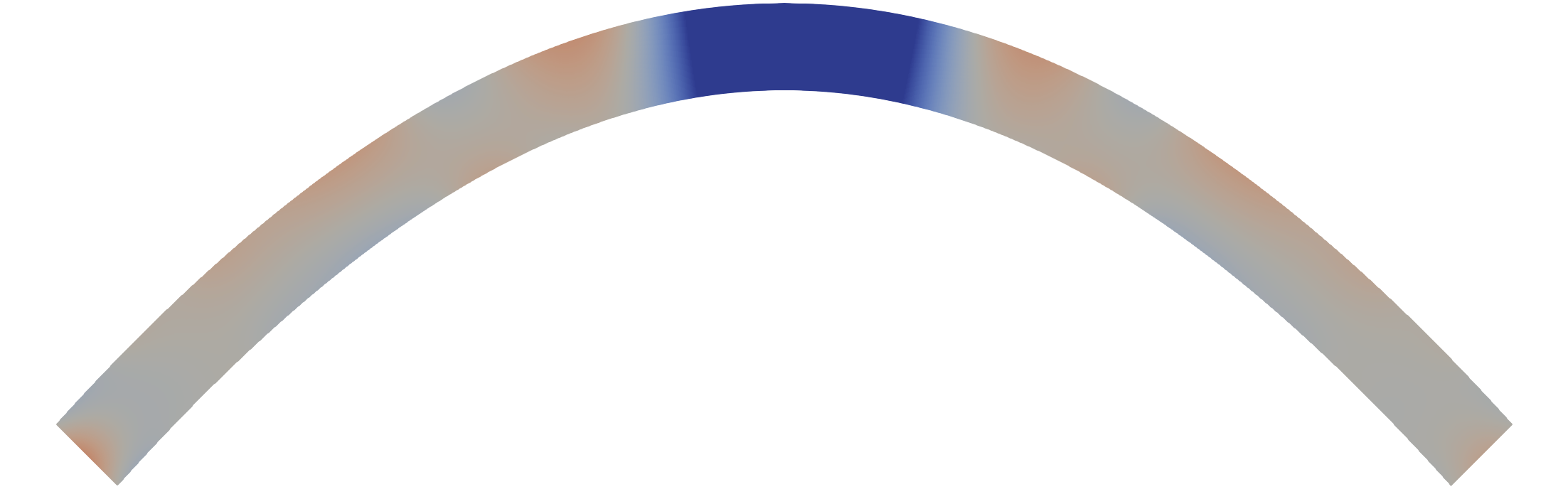}
\\[-11mm]
$t = 0.48$
\\[5mm]
\includegraphics[width=0.47\textwidth]{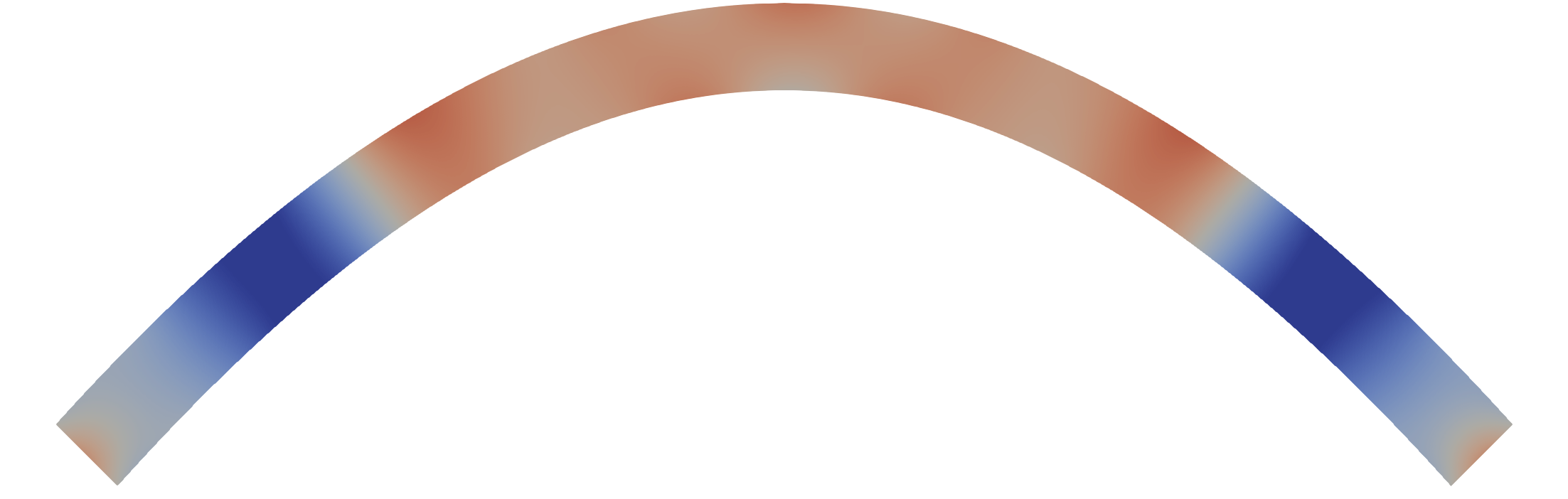}
\\[-11mm]
$t = 0.75$
\\[5mm]
\includegraphics[width=0.47\textwidth]{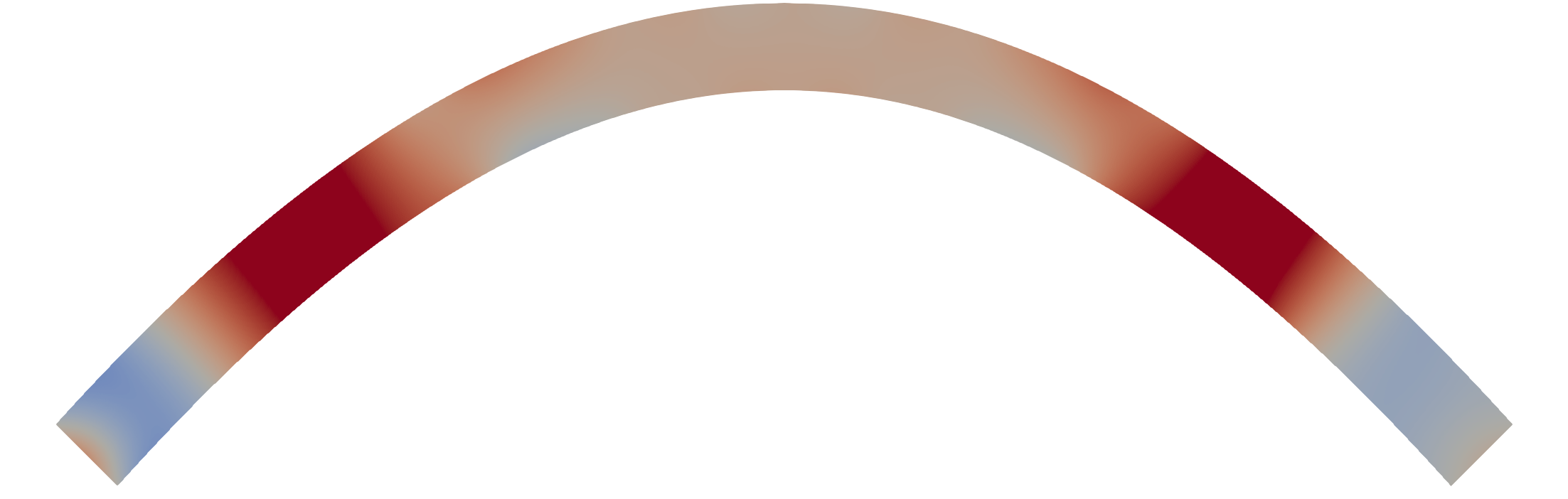}
\\[-11mm]
$t = 1$
\\[5mm]
\includegraphics[width=0.47\textwidth]{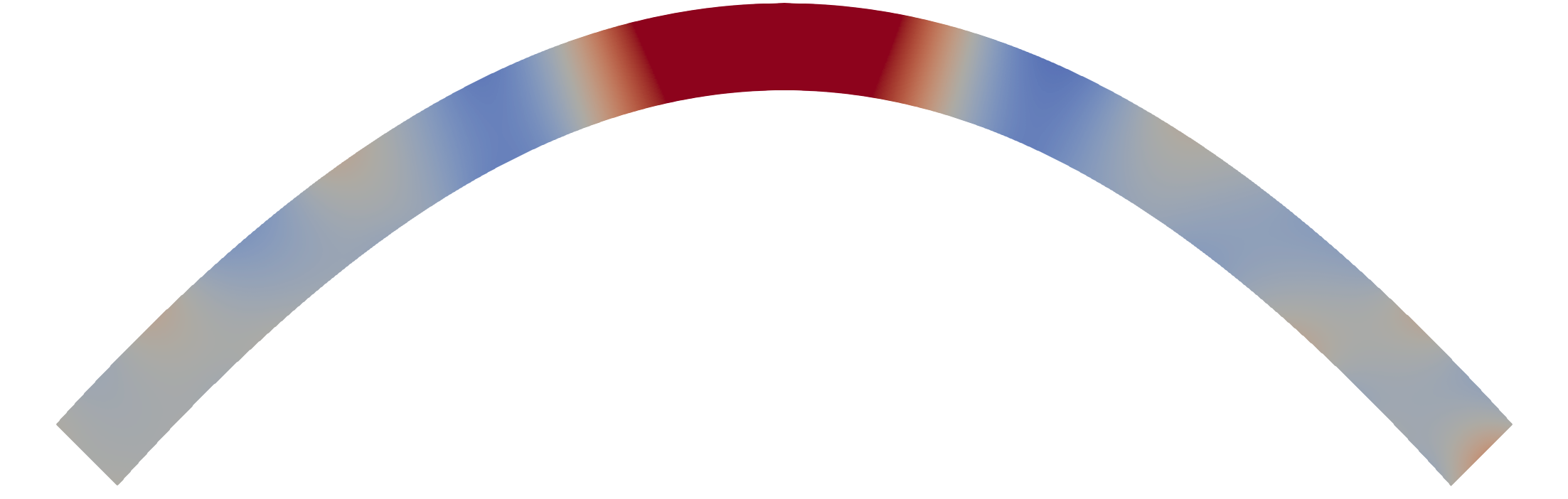}
\\[-11mm]
$t = 1.2$
\\[5mm]
\end{tabular}
\begin{tabular}{c}
\includegraphics[width=0.47\textwidth]{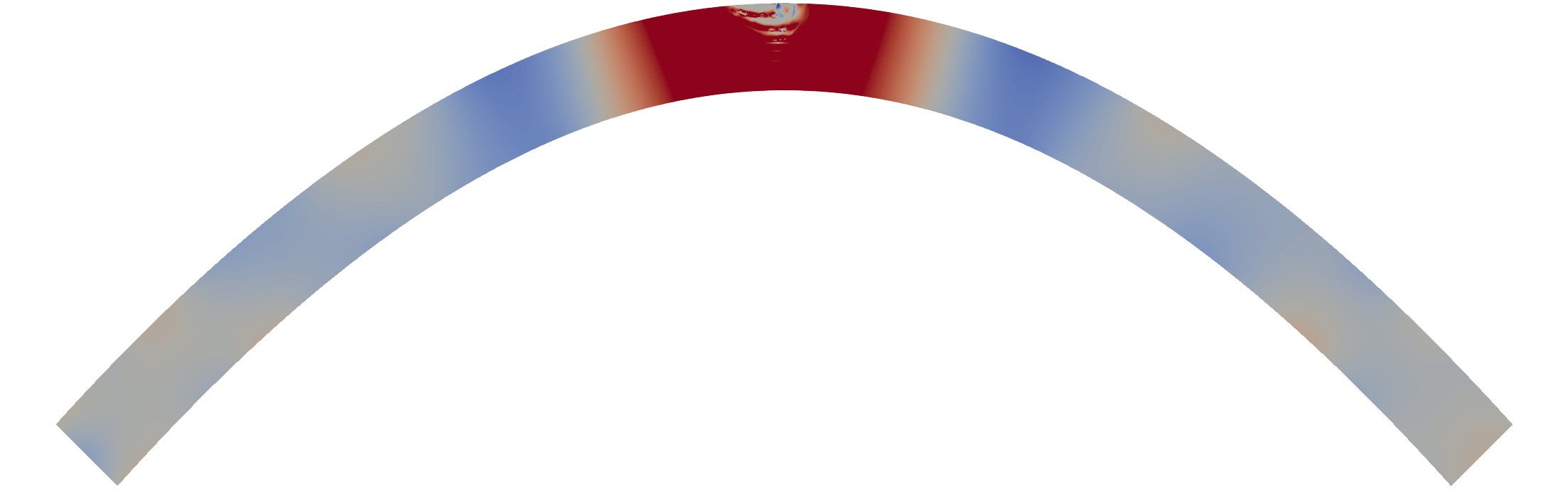}
\\[-11mm]
$t = 1.21$
\\[5mm]
\includegraphics[width=0.47\textwidth]{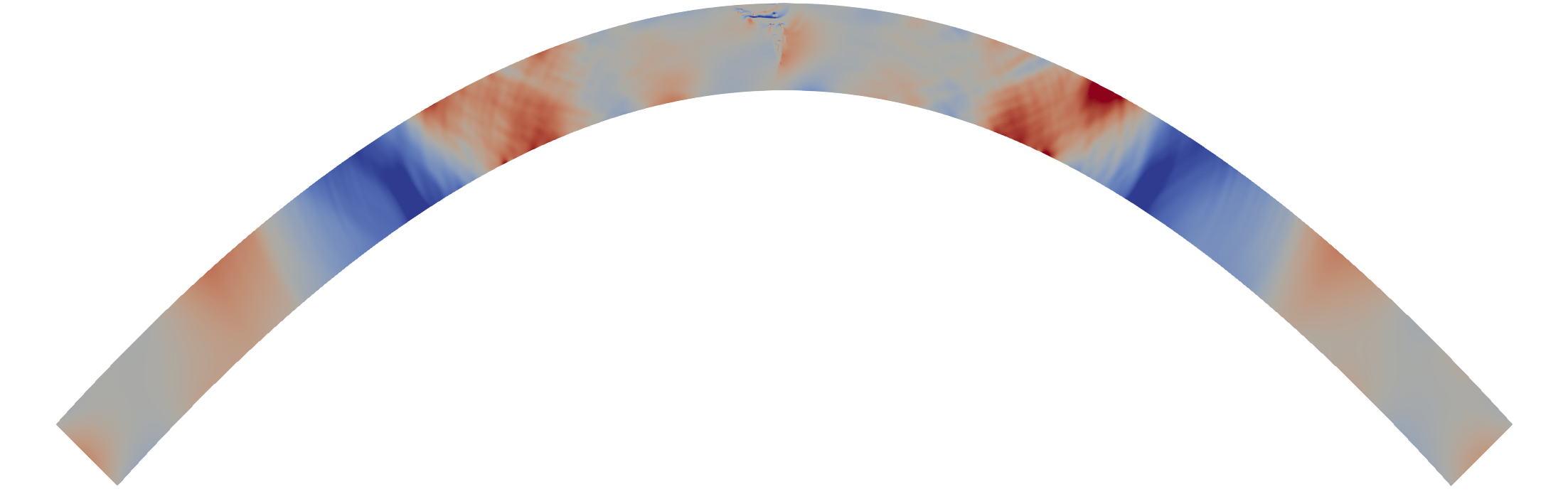}
\\[-11mm]
$t = 1.44$
\\[5mm]
\includegraphics[width=0.47\textwidth]{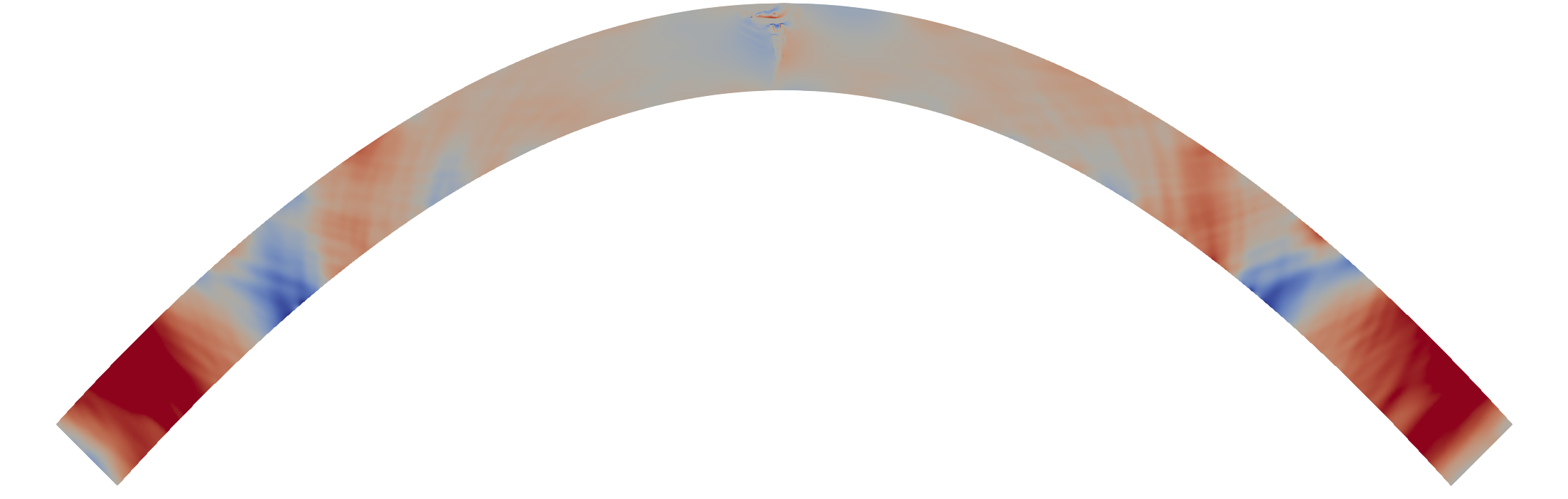}
\\[-11mm]
$t = 1.64$
\\[5mm]
\includegraphics[width=0.47\textwidth]{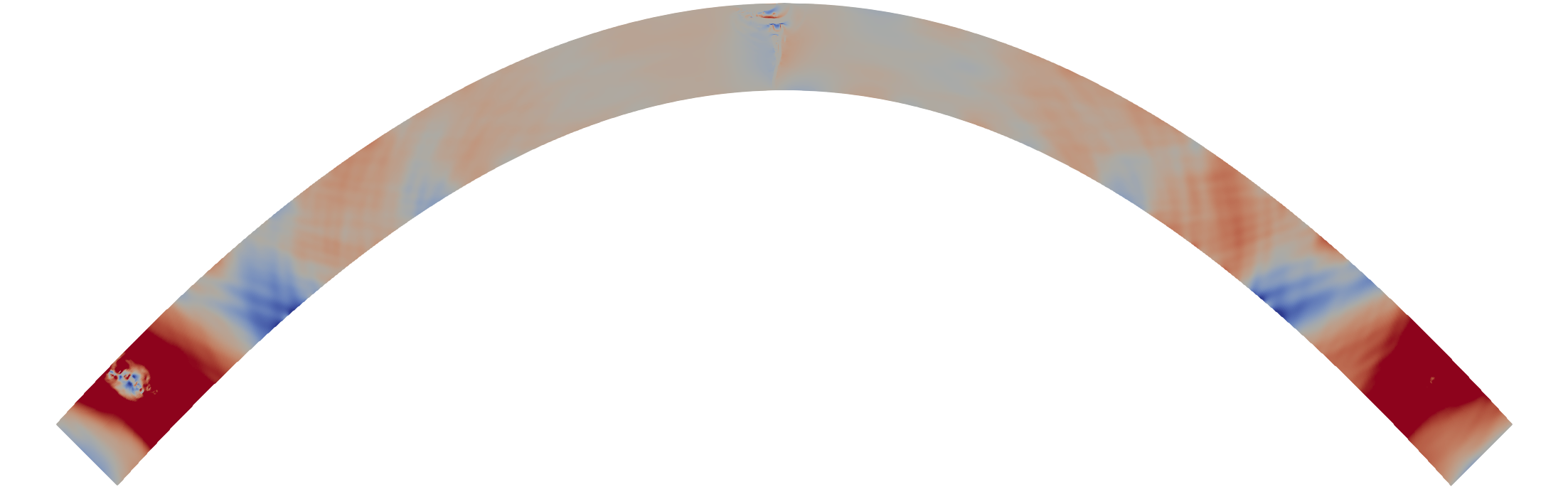}
\\[-11mm]
$t = 1.65$
\\[5mm]
\includegraphics[width=0.47\textwidth]{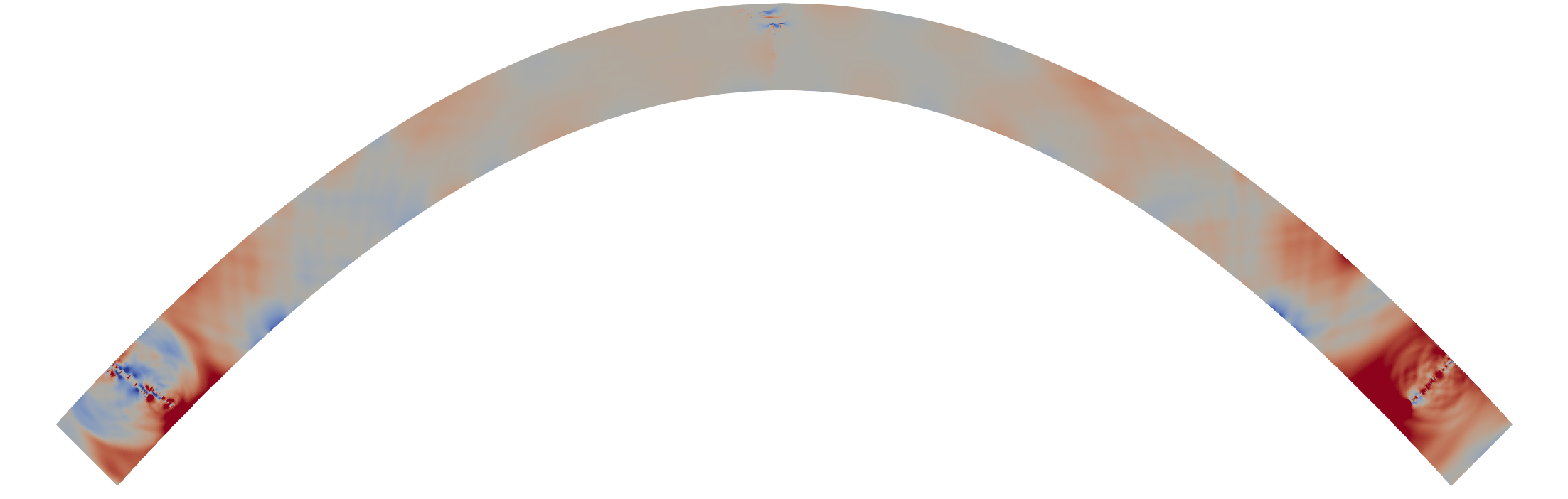}
\\[-11mm]
$t = 1.66$
\\[5mm]
\end{tabular}
\caption{Traveling pressure wave, $p=\operatorname{tr}\vec\sigma$ for $t\in [0,2]$, for the plane curved bar;\qquad\qquad\qquad  red -- compression, blue -- tension.}
\label{fig:trS}
\end{figure}

\gruen{In this numerical experiment, the following behavior is observed (illustrated by
snapshots of the solution in Fig.~\ref{fig:trS}): Compressional waves
are initiated for $t\in [0,{0.24}]$ by a pressure impulse from left and
right. Then, they travel, and for $t\in [0.75,1]$ the waves are reflected at the free
boundaries and result in tensional waves.  At $t=1.2$, by
superposition of the reflected waves, the tension increases, and the traction
forces become so large that the driving force $Y_\text{el}$ gets
positive and the fracture criterion is met; the phase field evolves
and the material breaks.  Since the domain is curved, the stress is
larger at the top side and the crack grows from top to bottom. For $t
> 1.2$, the waves are reflected at the crack interfaces and become
compressive.  Then, for $t\in [1.6,1.65]$ the waves are reflected
again at the free boundaries and turn into tension. Finally, for
$t=[1.65,1.71]$ by superposition of incoming and reflected waves, a
second fracture zone is initiated. Because of the non-symmetric setting, this secondary crack is slightly different at the left and the
right hand side.}

\clearpage


\begin{figure}[H]
\begin{tabular}{c}
\includegraphics[width=0.8\textwidth]{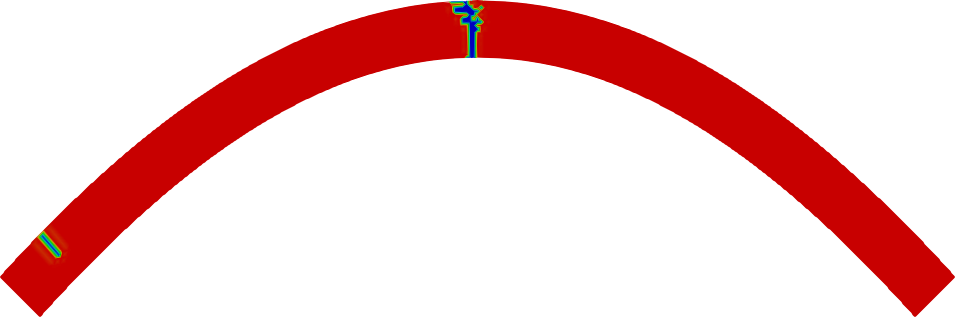} \\[-15mm]
$h = 2^{-8}$,
$|\mathcal K_h| = 4\,096$,\\
$\dim V^\text{cf} = 4\,369$, $\dim V_h^\text{dG} =  184\,320 $
\\[5mm]
\includegraphics[width=0.8\textwidth]{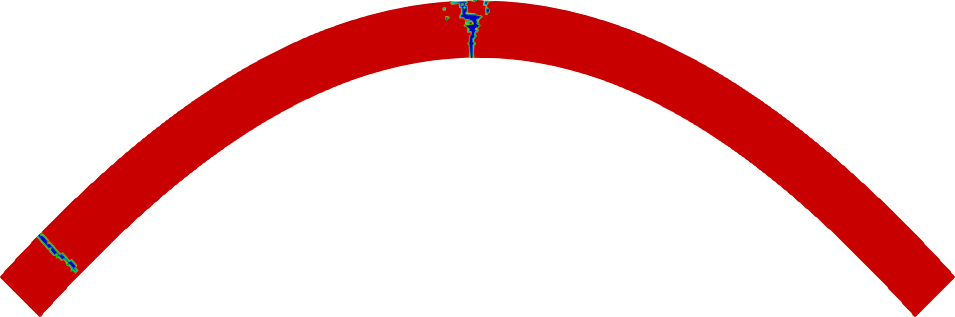} \\[-15mm]
$h = 2^{-9}$,  $|\mathcal K_h| = 16\,384$,\\
$\dim V^\text{cf} = 16\,929$, $\dim V_h^\text{dG} =737\,280 $
\\[5mm]
\includegraphics[width=0.8\textwidth]{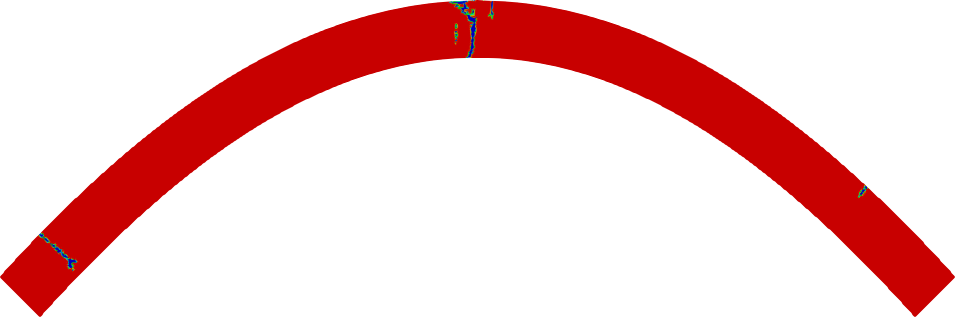} \\[-15mm]
$h = 2^{-10}$, $|\mathcal K_h| = 65\,536$,\\
$\dim V^\text{cf} = 66\,625$, $\dim V_h^\text{dG} =1\,310\,720 $
\\[5mm]
\includegraphics[width=0.8\textwidth]{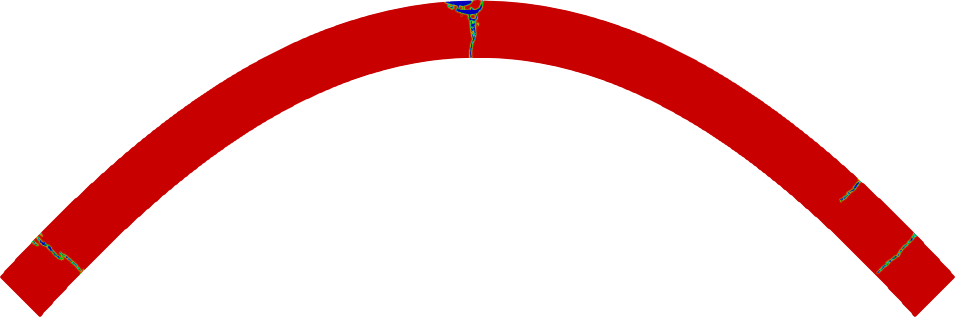} \\[-15mm]
$h = 2^{-11}$,
$|\mathcal K_h| = 262\,144$,\\
$\dim V^\text{cf} = 264\,321$,
$\dim V_h^\text{dG} = 5\,242\,880$
\\[5mm]
\end{tabular}
\caption{Phase-field fracture approximation   $s_{h,\text{inf}}$ at $t=2$ computed with time step sizes
  ${\vartriangle} t_\text{el}= 0.001$ and ${\vartriangle} t_\text{pf}=0.0005$
  on different meshes with mesh size $h$,
  number of elements $|\mathcal K_h|$, dimension
of the finite element spaces for the phase field $\dim V^\text{cf}$ and the wave system
$\dim V_h^\text{dG}$.
}
\label{fig:pf}
\end{figure}

\gruen{The resulting
phase-field approximation of the fracture evolution is shown in
Fig.~\ref{fig:pf}. We use a fixed length scale of $l_\text{c}= 0.0005$
which determines the fracture regularization by the phase-field
approach. This requires to use a mesh size of $h= 2^{-8}$ to
resolve the main features of the fracture pattern including primary
and the secondary cracks.  Then, using finer meshes, we observe than
the overall setting is converging, even though the cracks are more
complex on finer meshes. The accuracy of the
computation is limited by the a priori choice of the phase-field
regularization given by the length scale parameter, resulting in a
diffusive approximation of the crack. Moreover, we observe
that larger values of the scaling parameter $M_\text{geom}$ or of the
length scale $l_\text{c}$ increase the dissipation
and so, after the primary crack,
the remaining mechanical energy is not sufficient for the
initiation of further cracks.  }

\gruen{
In this numerical experiment we use on level $m=8,9,10,11$ the mesh
size $h= 2^{-m}$ with $2^{m+4}$ quadrilateral elements and 45 degrees of
freedom per element for the velocity and stress
approximation. The smallest computation on level $m=8$ with 2\,000
time steps runs on a laptop within 14 hours,
so that we need about 25 seconds per time step.
The finer computations were realized on a parallel computing cluster
with different numbers of cores.
}

\begin{figure}[H]
  \centering\small
  \begin{tabular}{p{0.1\textwidth}p{0.1\textwidth}p{0.1\textwidth}p{0.1\textwidth}p{0.1\textwidth}p{0.1\textwidth}p{0.1\textwidth}p{0.1\textwidth}}
    \includegraphics[width=0.12\textwidth]{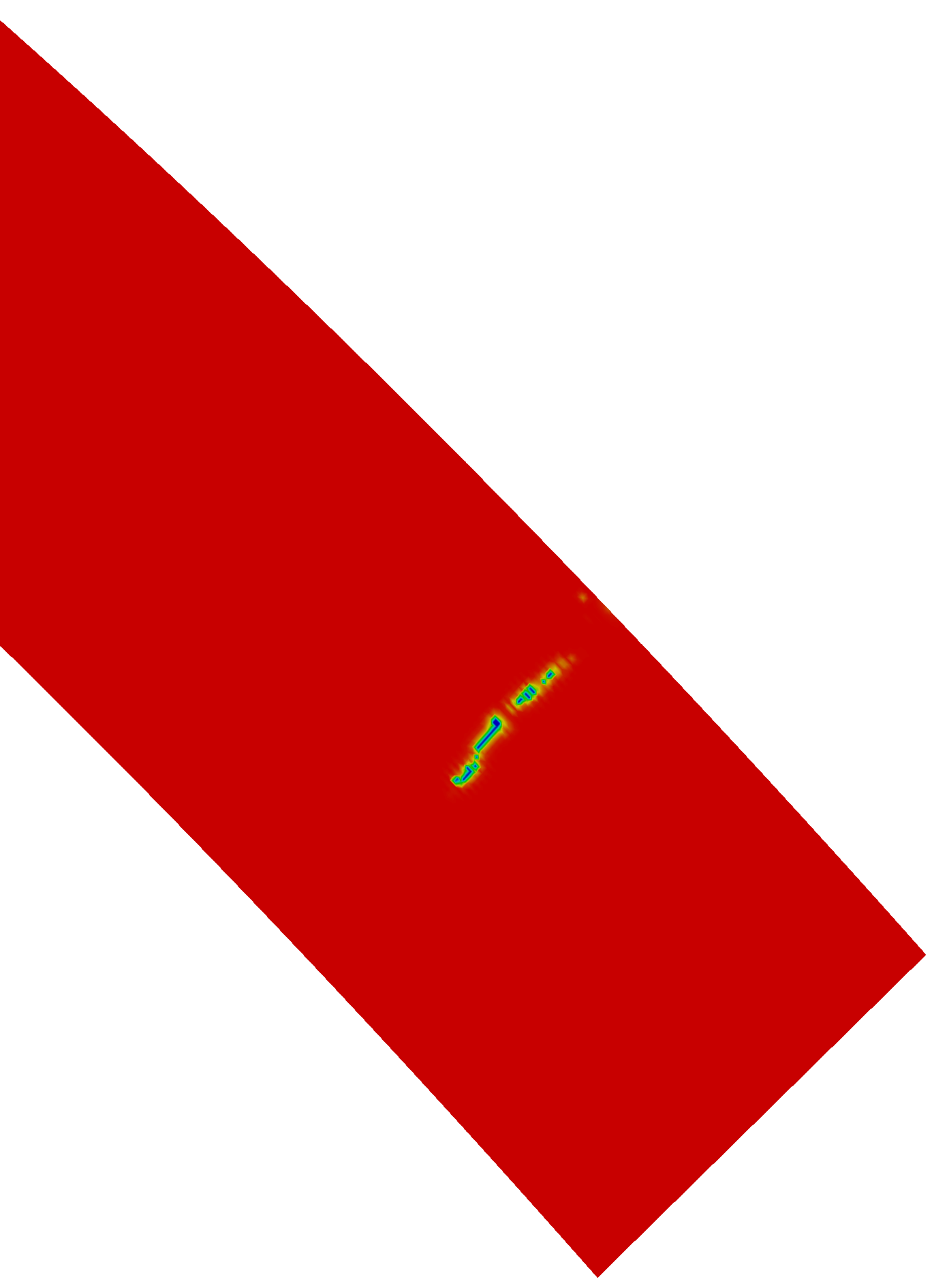}
    &\includegraphics[width=0.12\textwidth]{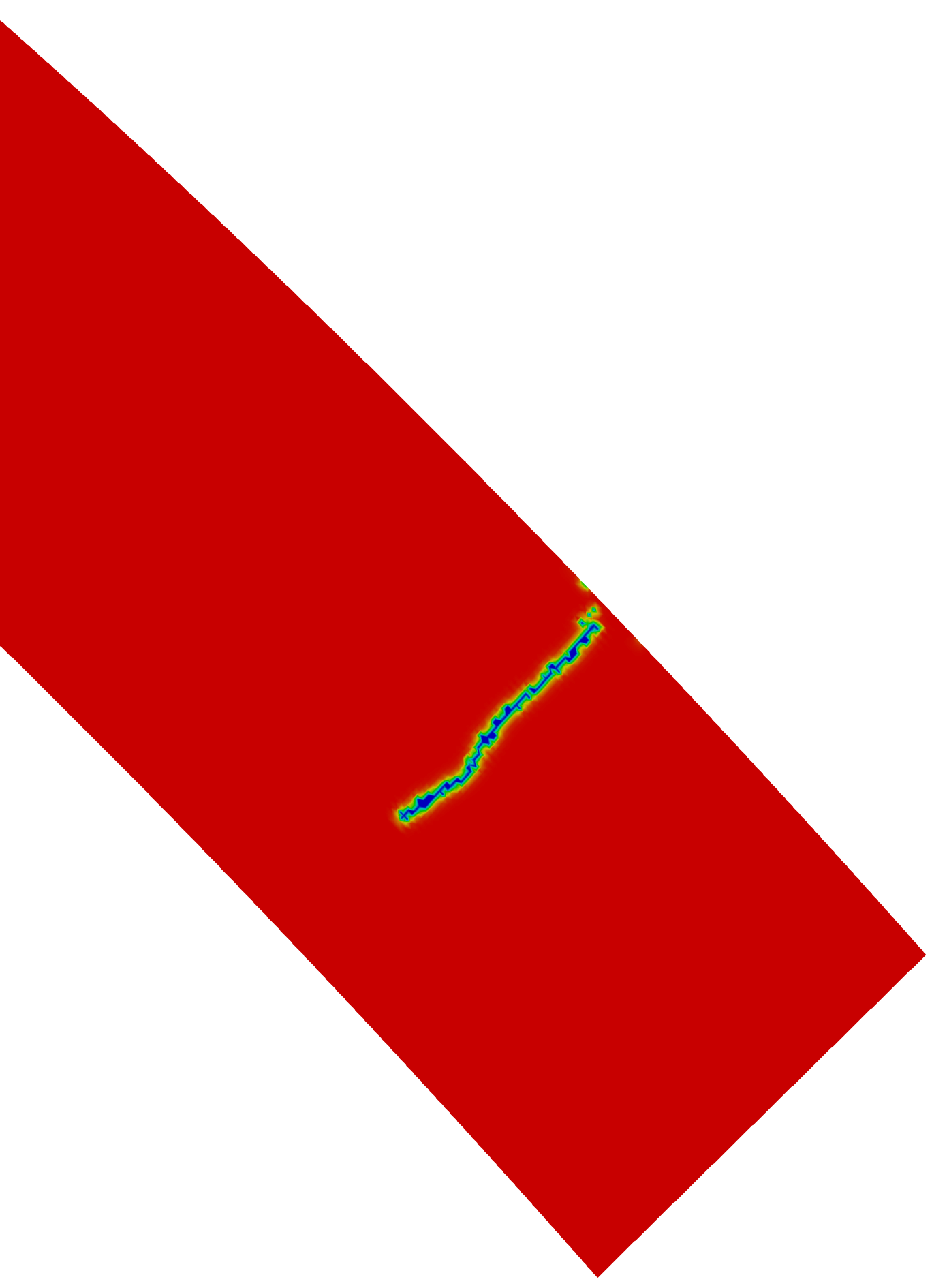}
    &\includegraphics[width=0.12\textwidth]{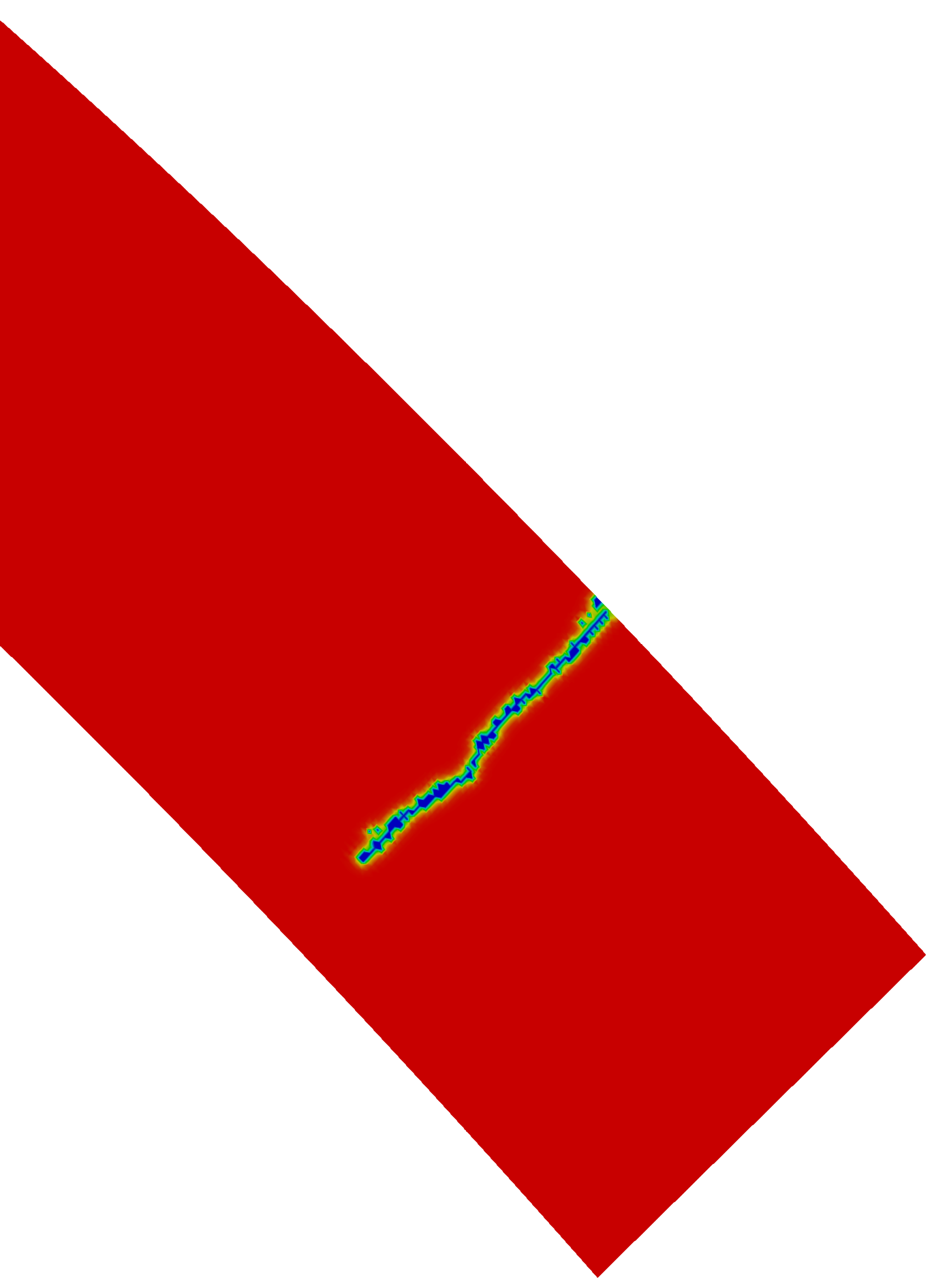}
    &\includegraphics[width=0.12\textwidth]{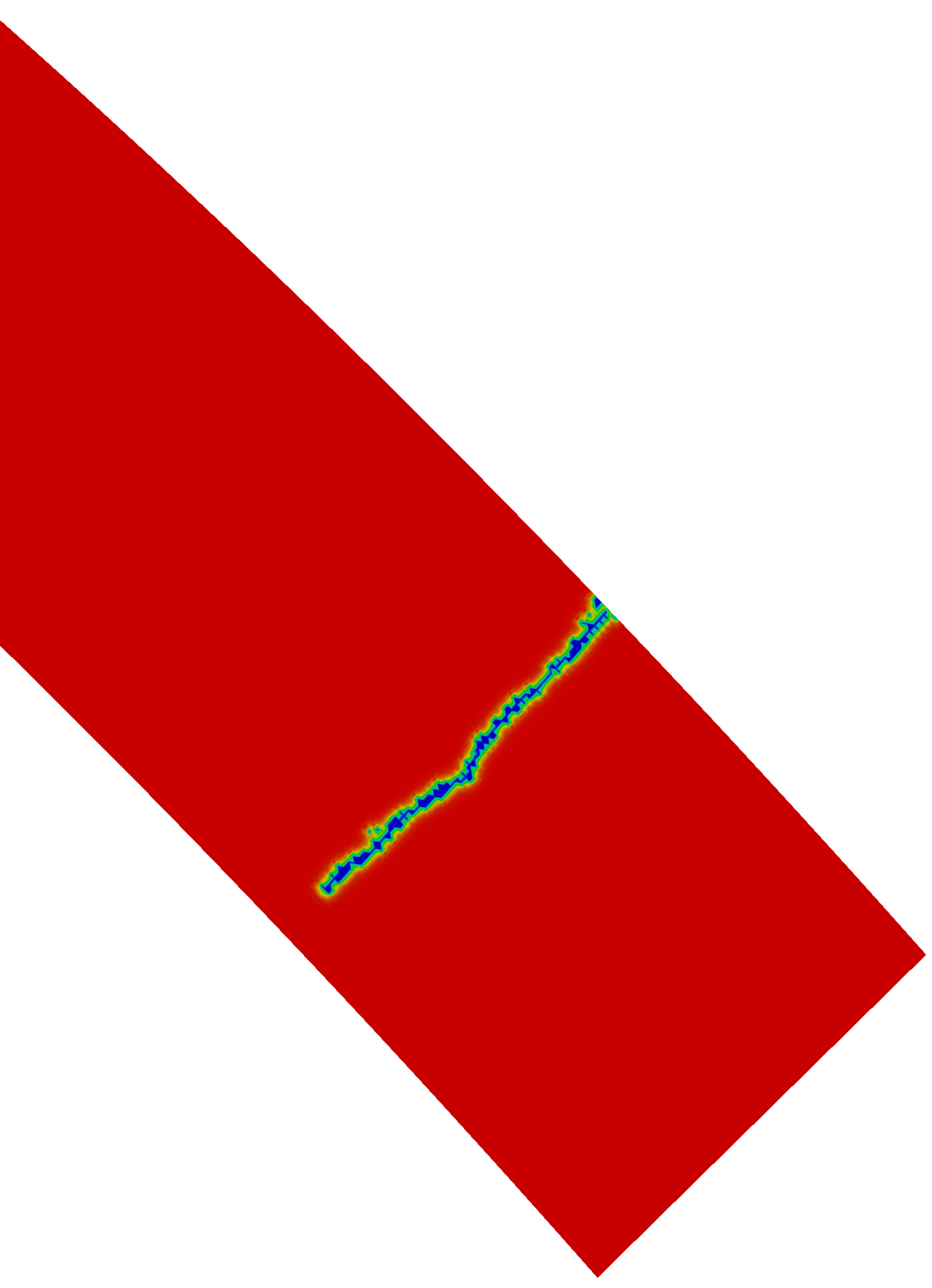}
    &\includegraphics[width=0.12\textwidth]{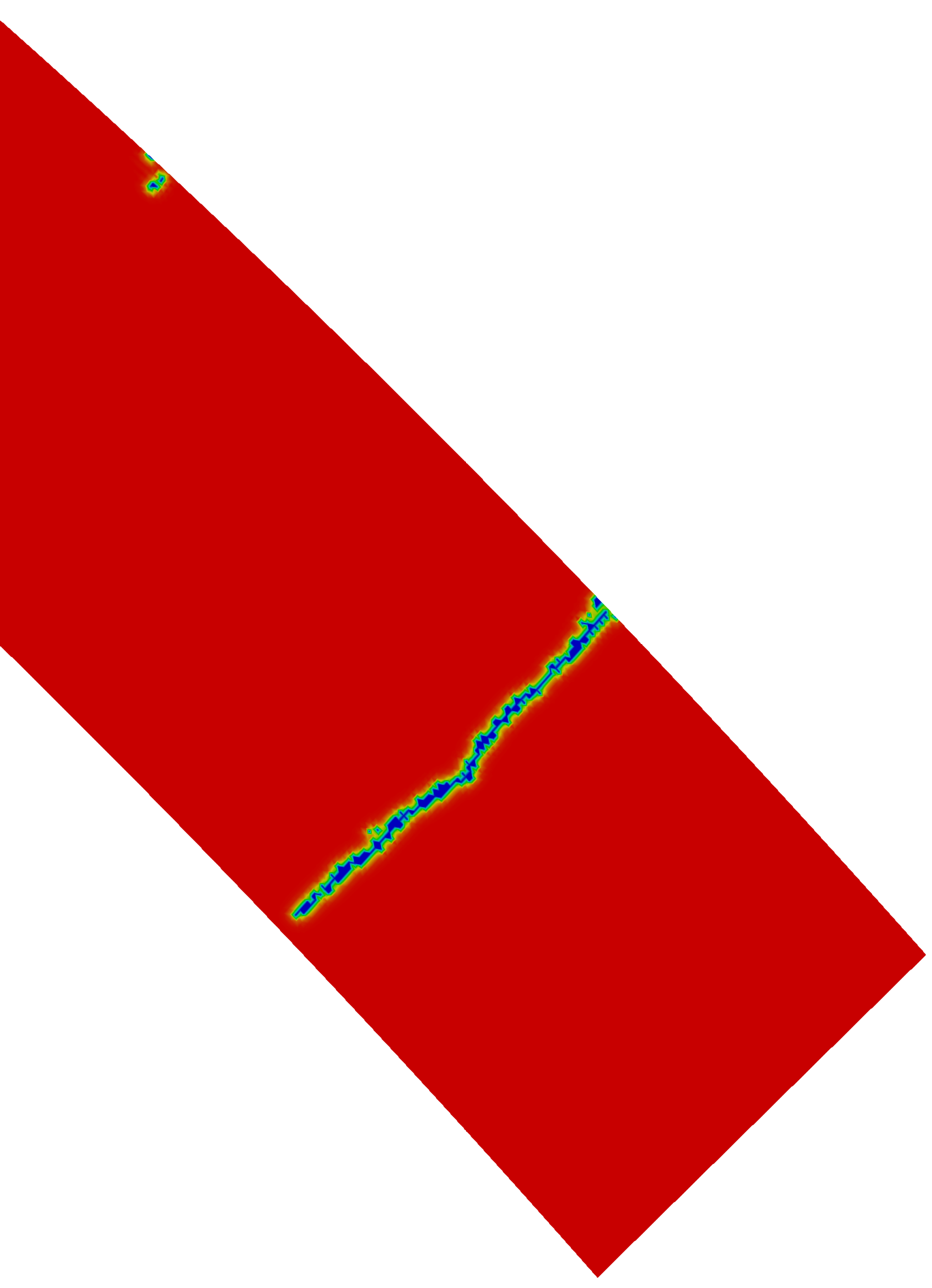}
    &\includegraphics[width=0.12\textwidth]{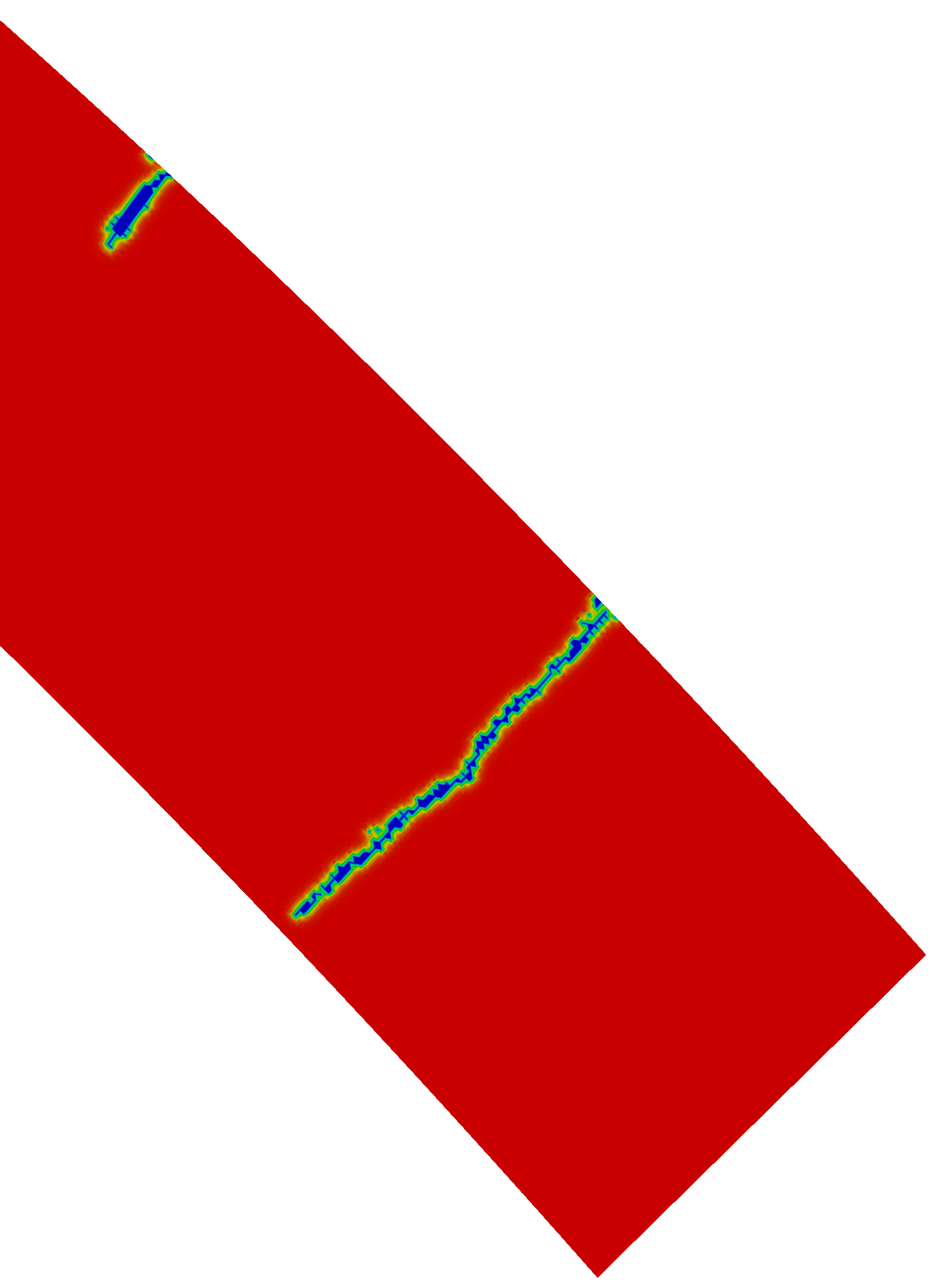}
    &\includegraphics[width=0.12\textwidth]{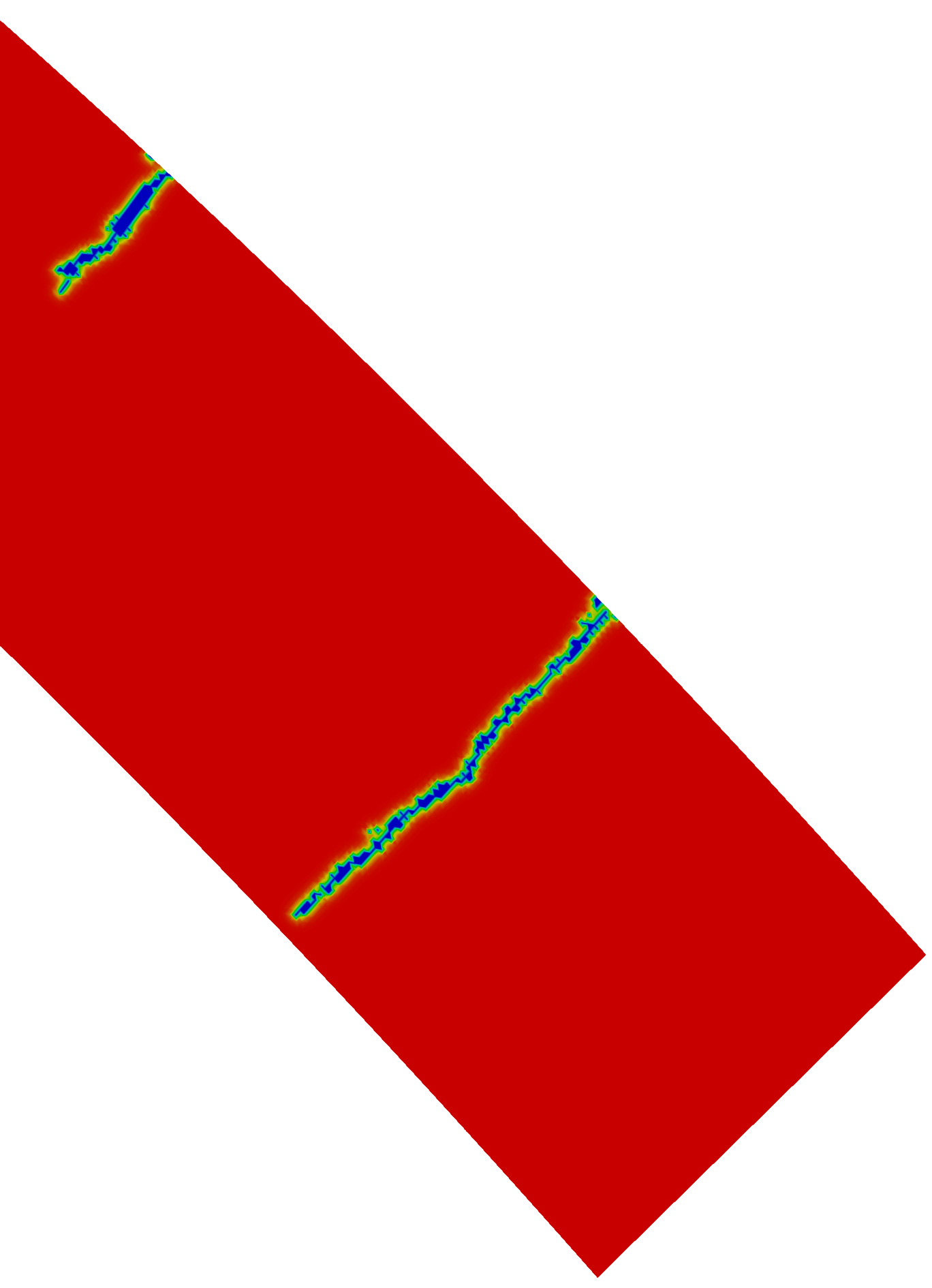}
    &\includegraphics[width=0.12\textwidth]{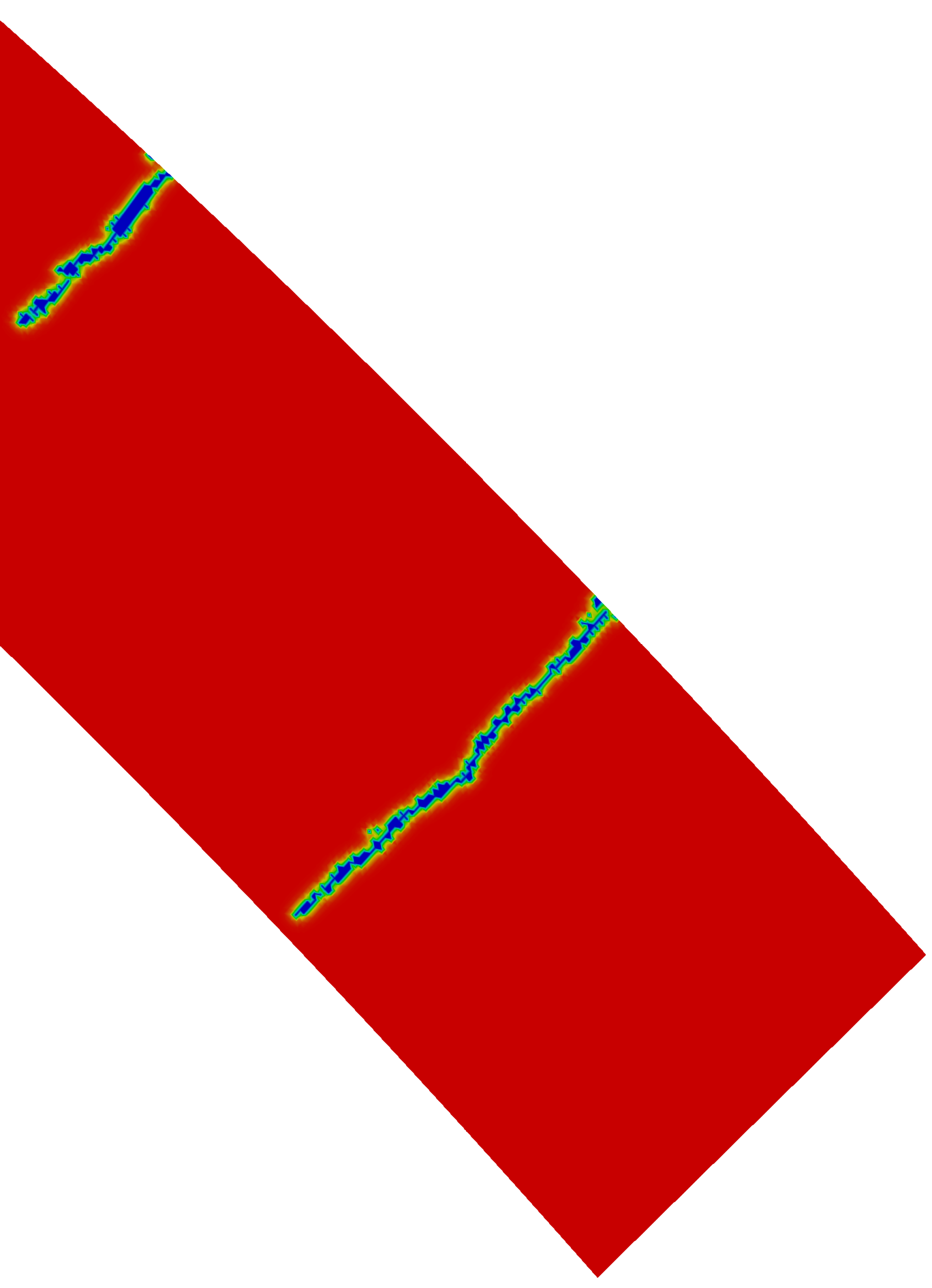}
    \\[-1mm]
    $\quad t= 1.660$&
    $\quad t= 1.670$&
    $\quad t= 1.680$&
    $\quad t= 1.690$&
    $\quad t= 1.700$&
    $\quad t= 1.703$&
    $\quad t= 1.708$&
    $\quad t= 1.710$
    \end{tabular}
  \centering
  \caption{Evolution of a secondary crack close to the right boundary, $s_{h,\text{inf}}$ for $t\in [1.66,1.71]$; the crack tip propagates with approximately 50\% of the wave speed.}
\label{fig:cracks}
\end{figure}

In Fig.~\ref{fig:cracks} the evolution of the secondary crack at the
right hand side of the model is displayed \gruen{on the finest mesh}. Here it can be
observed that the crack propagates with approximately 50\% of the wave
speed.
Such cracks, with a velocity of about the Rayleigh-wave speed at the
crack tip, are typical for dynamic brittle fracture,
cf.~\cite[Chap.~7.4]{freund1998dynamic}.


\clearpage

\clearpage

\subsection{A 3D experiment}
\begin{figure}[H]
  \begin{tabular}{cc}
    \\[-3mm]
    $t_n=0.125$, $n = 5\,000$
    \\[-10mm]
    &
    \includegraphics[width=0.5\textwidth]{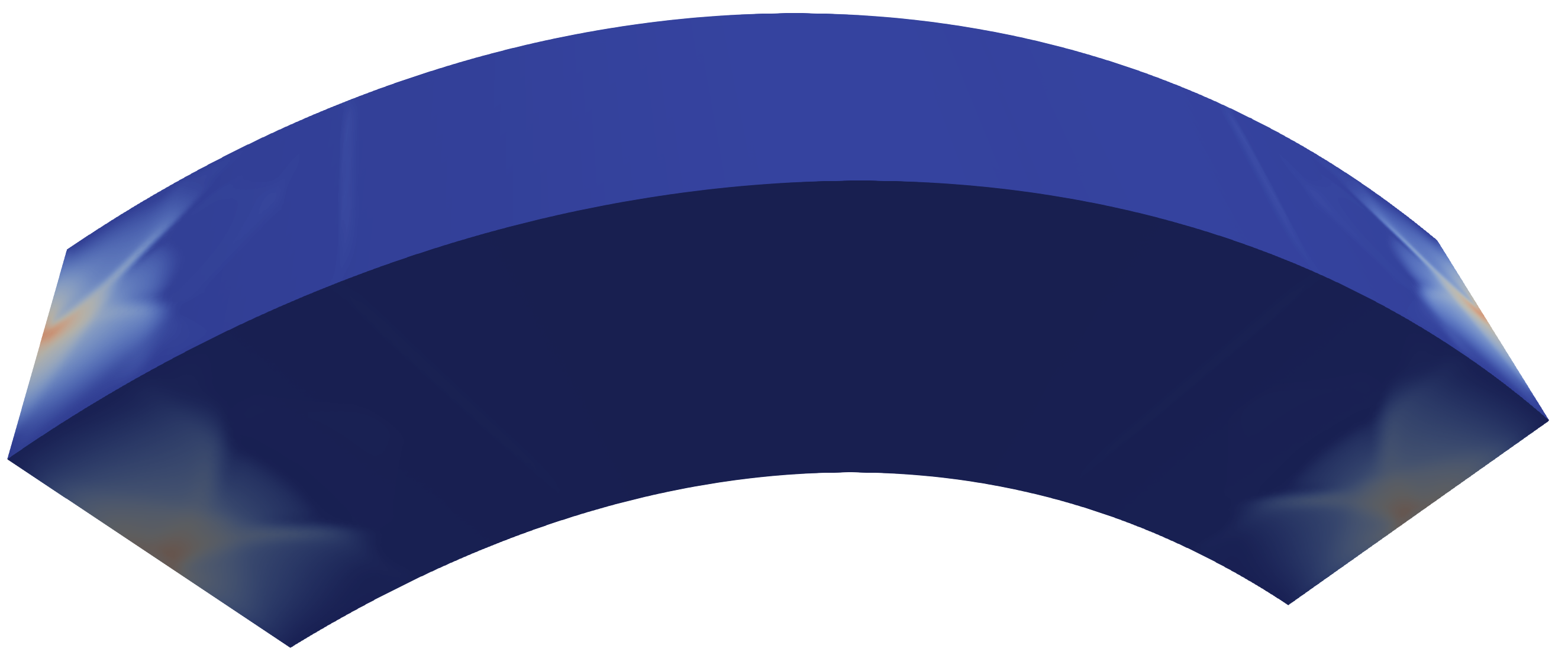}
    \\[3.5mm]
    $t_n=0.25$, $n = 10\,000$
    \\[-10mm]
    &
    \includegraphics[width=0.5\textwidth]{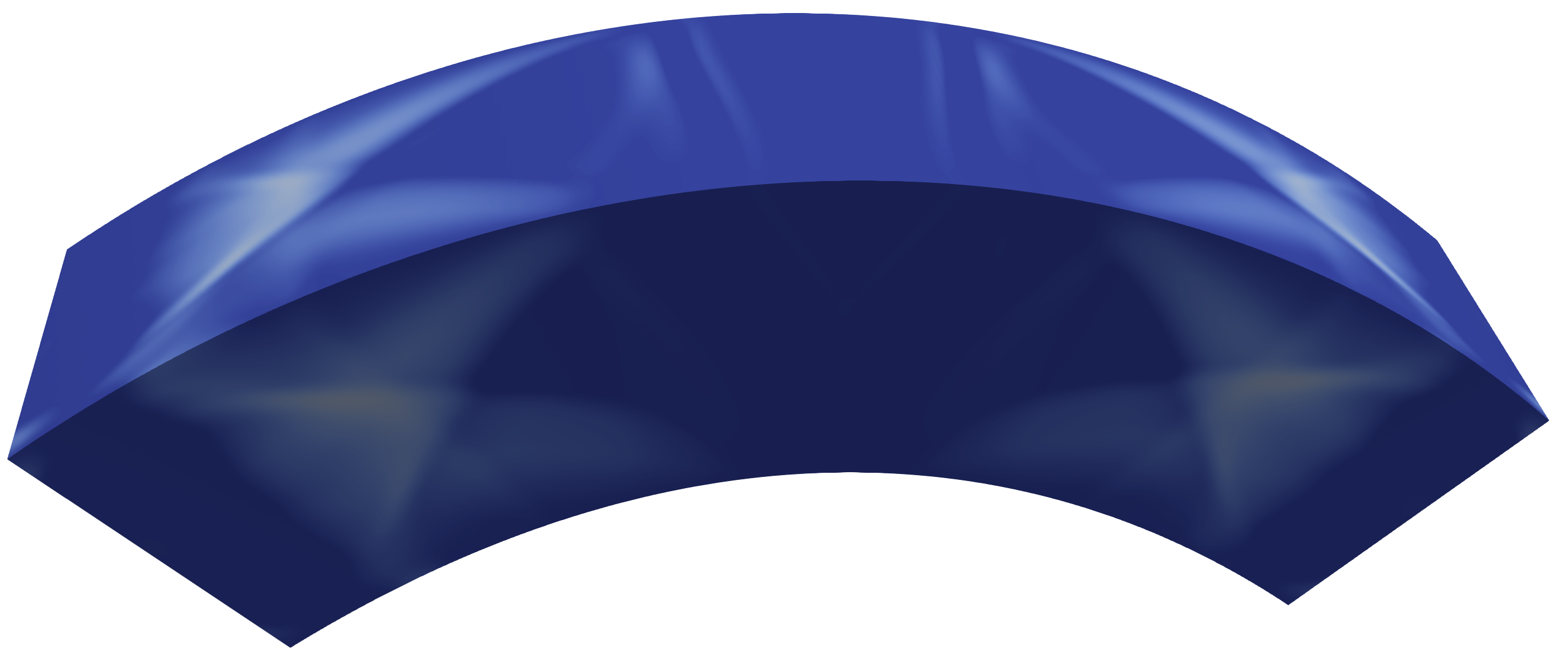}
    \\[3.5mm]
    $t_n=0.375$, $n = 15\,000$
    \\[-10mm]
    &
    \includegraphics[width=0.5\textwidth]{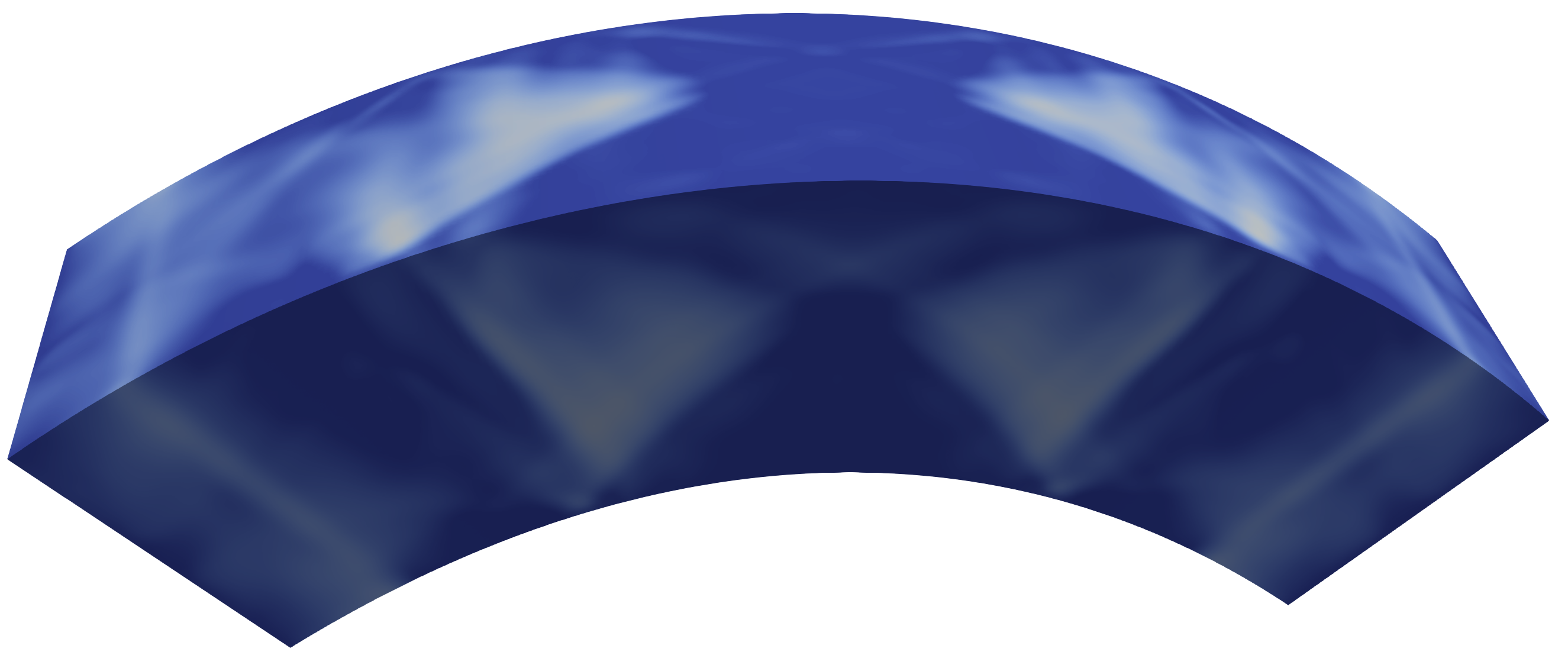}
    \\[3.5mm]
    $t_n=0.5$, $n = 20\,000$
    \\[-10mm]
    &
    \includegraphics[width=0.5\textwidth]{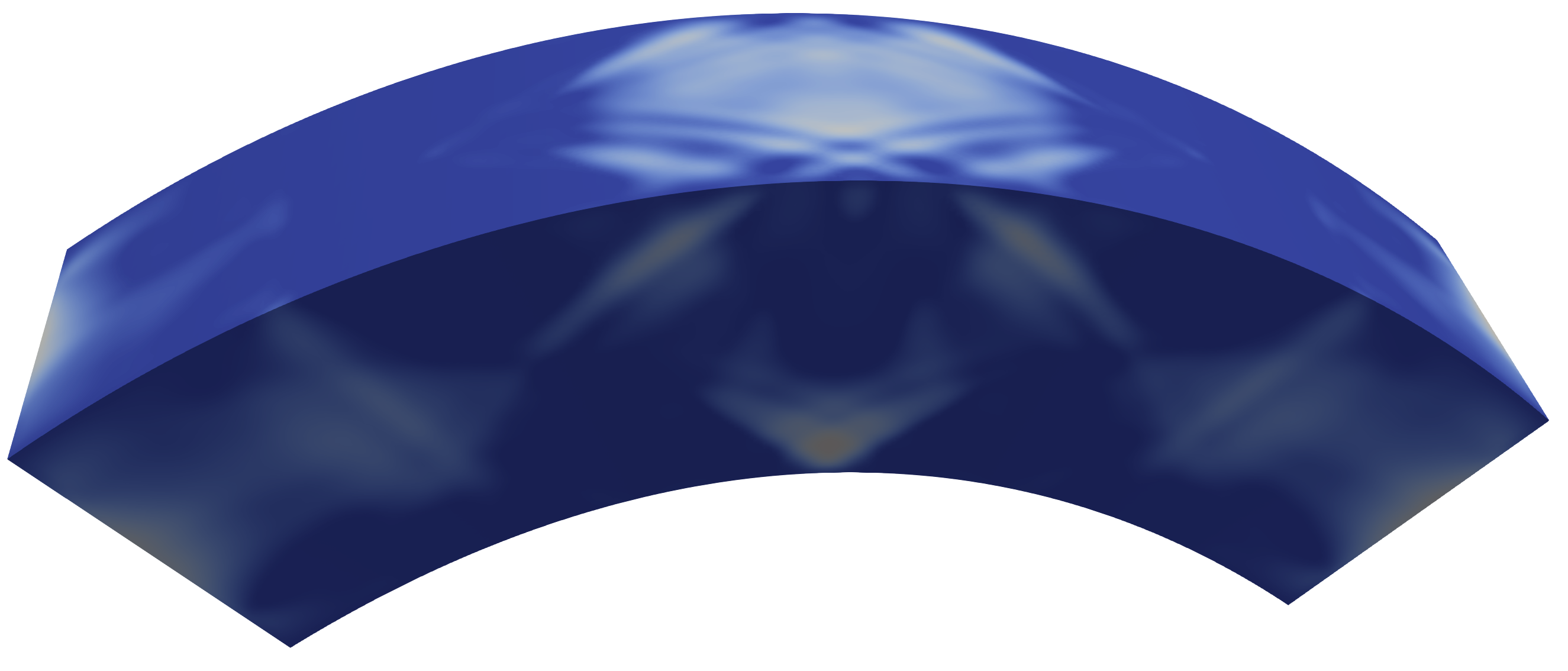}
    \\[3.5mm]
    $t_n=0.625$, $n=25\,000$
    \\[-10mm]
    &
    \includegraphics[width=0.5\textwidth]{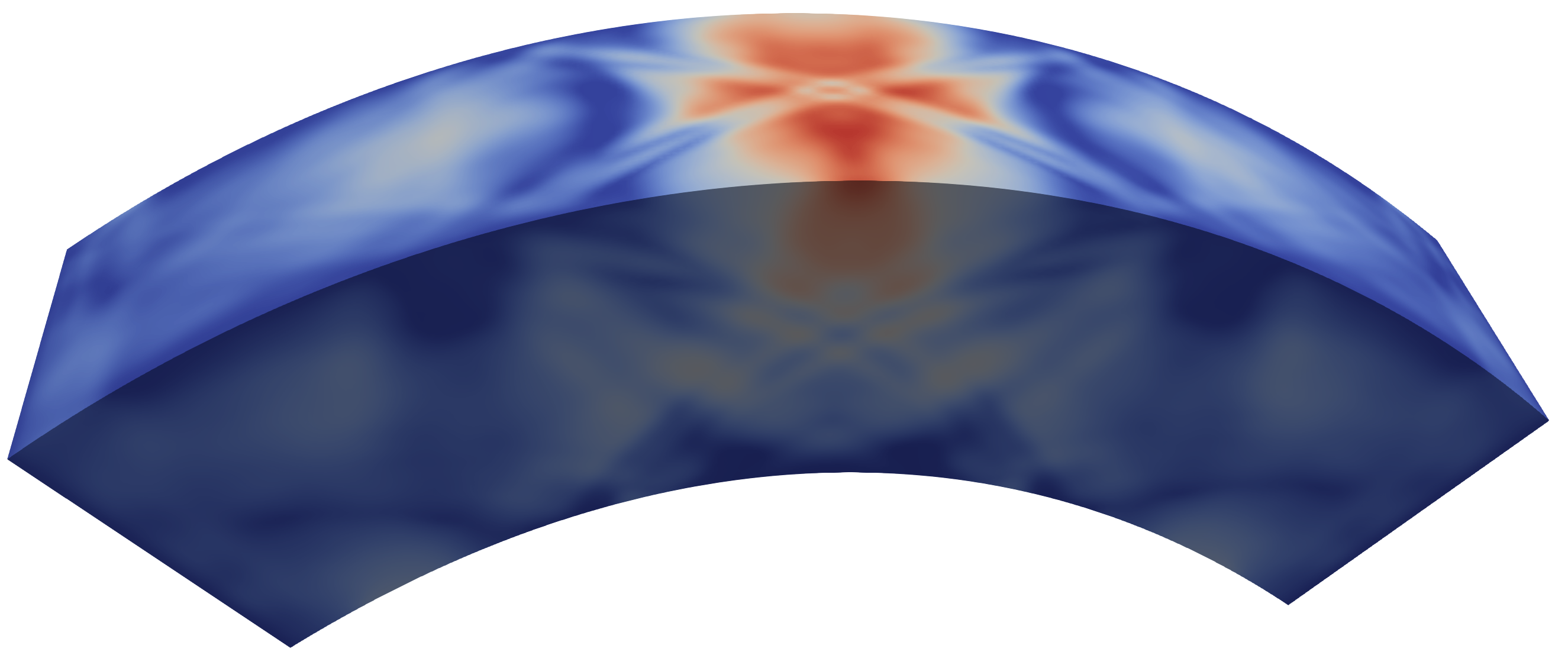}
    \\[3.5mm]
    $t_n=0.75$, $n=30\,000$
    \\[-10mm]
    &
    \includegraphics[width=0.5\textwidth]{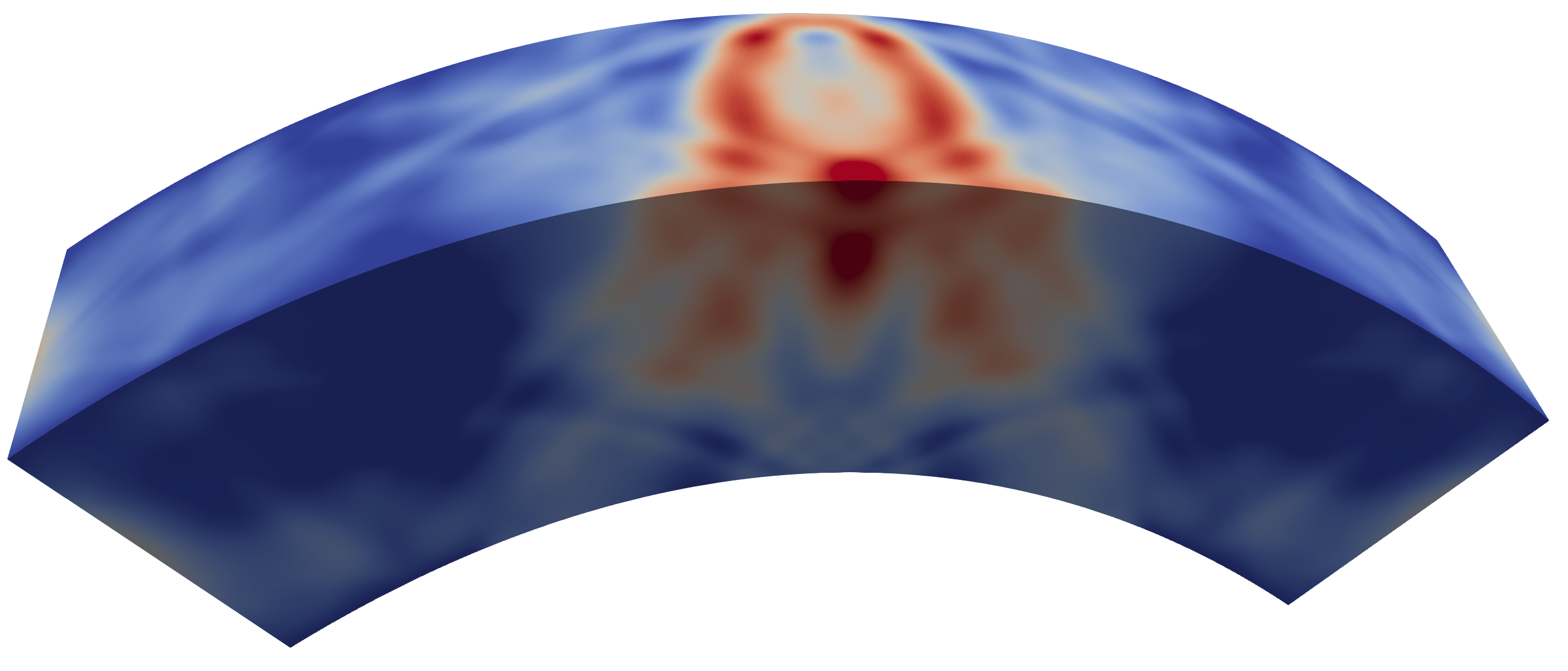}
    \\[-3.5mm]
\end{tabular}
\caption{Principle stress evolution $\sigma_\text{I}$ for the 3D configuration
on 1\,048\,576 hexahedra.
}
\label{fig:3dwave}
\end{figure}

\enlargethispage{0.5cm}
\gruen{
Here we show that the observed qualitative behavior transfers to three
space dimensions. The configuration is now somewhat simplified, i.e., the
curved bar-like domain is shorter, and the pressure impulse is
modified so that already the first superposition of the waves
generates tension. At the front side of the bar, the crack criterion
is met first.  The wave propagation and the evolution of cracks
patterns at the surface are illustrated in Fig.~\ref{fig:3dwave}.  }

\clearpage

\gruen{
In this simulation we use a fixed time step size
${\vartriangle} t_\text{el}= {\vartriangle} t_\text{pf}
= 0.000025$, a trilinear DG finite element space
for the elastic system with 75\,497\,472 degrees of freedom,
and a phase-field approximation with 1\,085\,825 degrees of freedom
on a hexahedral mesh.
The full simulation requires approx.~20 hours
on 2048 parallel computing cores on HoreKA
(\href{https://www.nhr.kit.edu/userdocs/horeka}{https://www.nhr.kit.edu/userdocs/horeka}),
so that per time step we need only a few seconds for the
parallel preconditioned GMRES solver.
}

\begin{figure}[H]
  \begin{tabular}{cc}
    \\[1mm]
    $t_n=0.75$, $n=30\,000$
    \\[-10mm]
    &
    \includegraphics[width=0.7\textwidth]{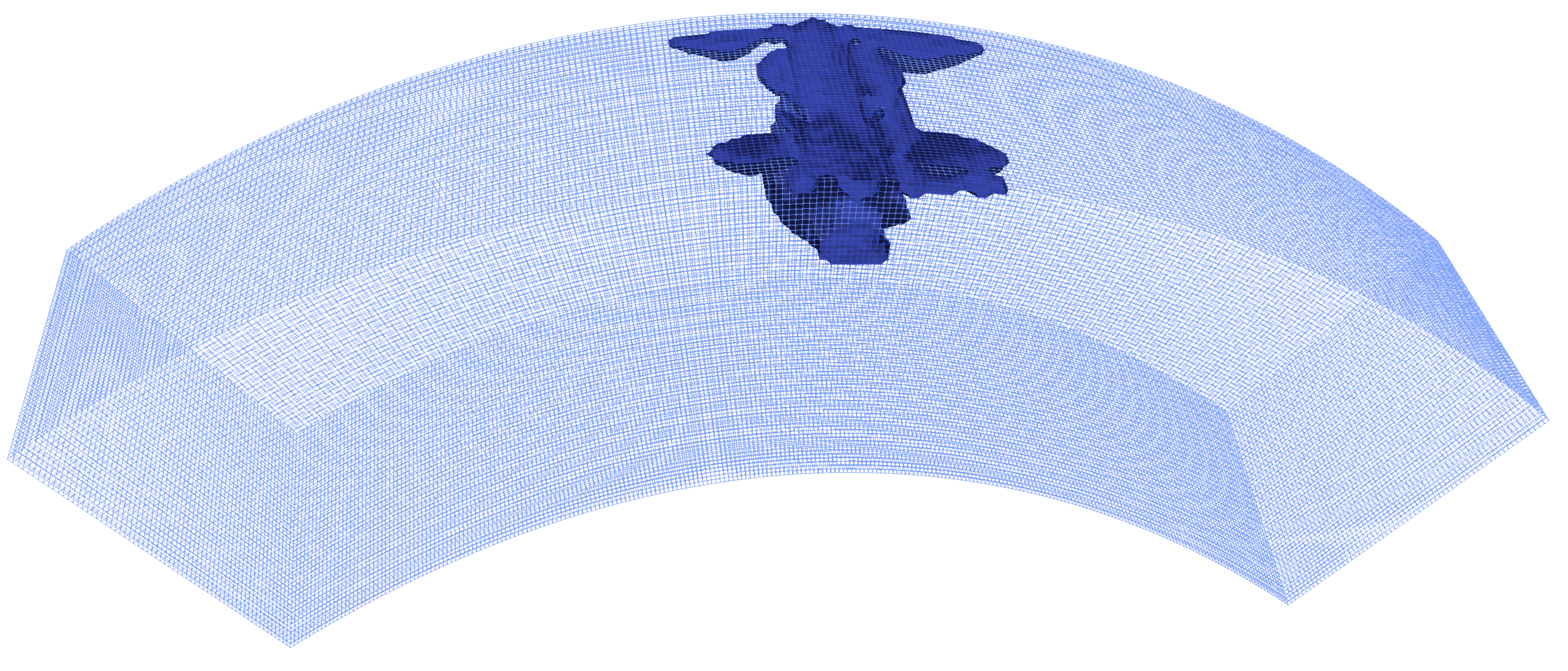}
    \\[-1mm]
\end{tabular}
  \caption{Phase-field approximation $s_{h,\text{inf}}$ at final time $t=0.75$.}
\label{fig:3dcrack}
\end{figure}

\begin{figure}[H]
  \begin{tabular}{ccccc}
    \multicolumn{5}{c}{phase-field approximation 
      with  1\,085\,825 DoFs on  1\,048\,576 hexahedra}
    \\[2mm]
    \includegraphics[width=0.18\textwidth]{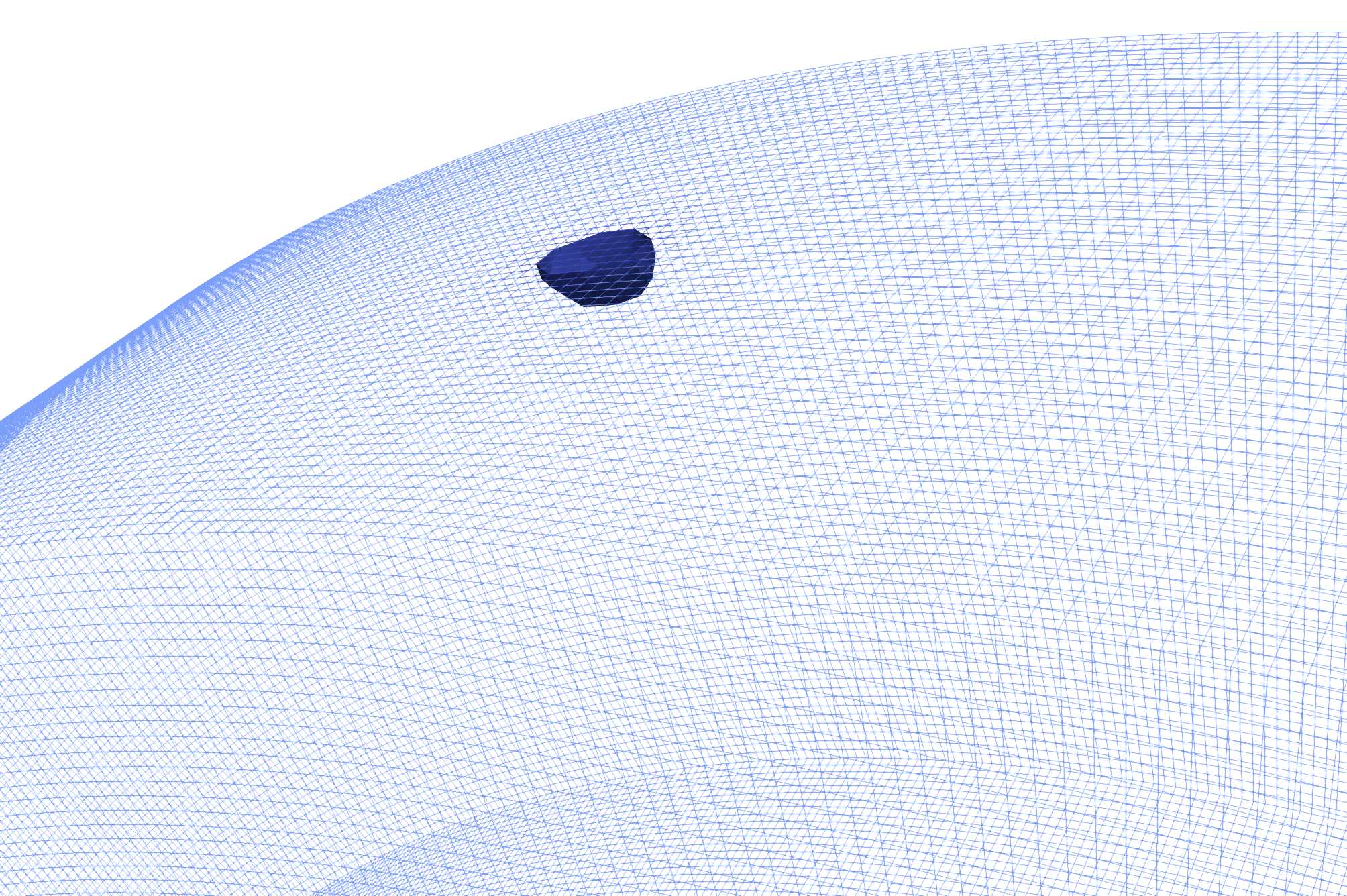}
    &\includegraphics[width=0.18\textwidth]{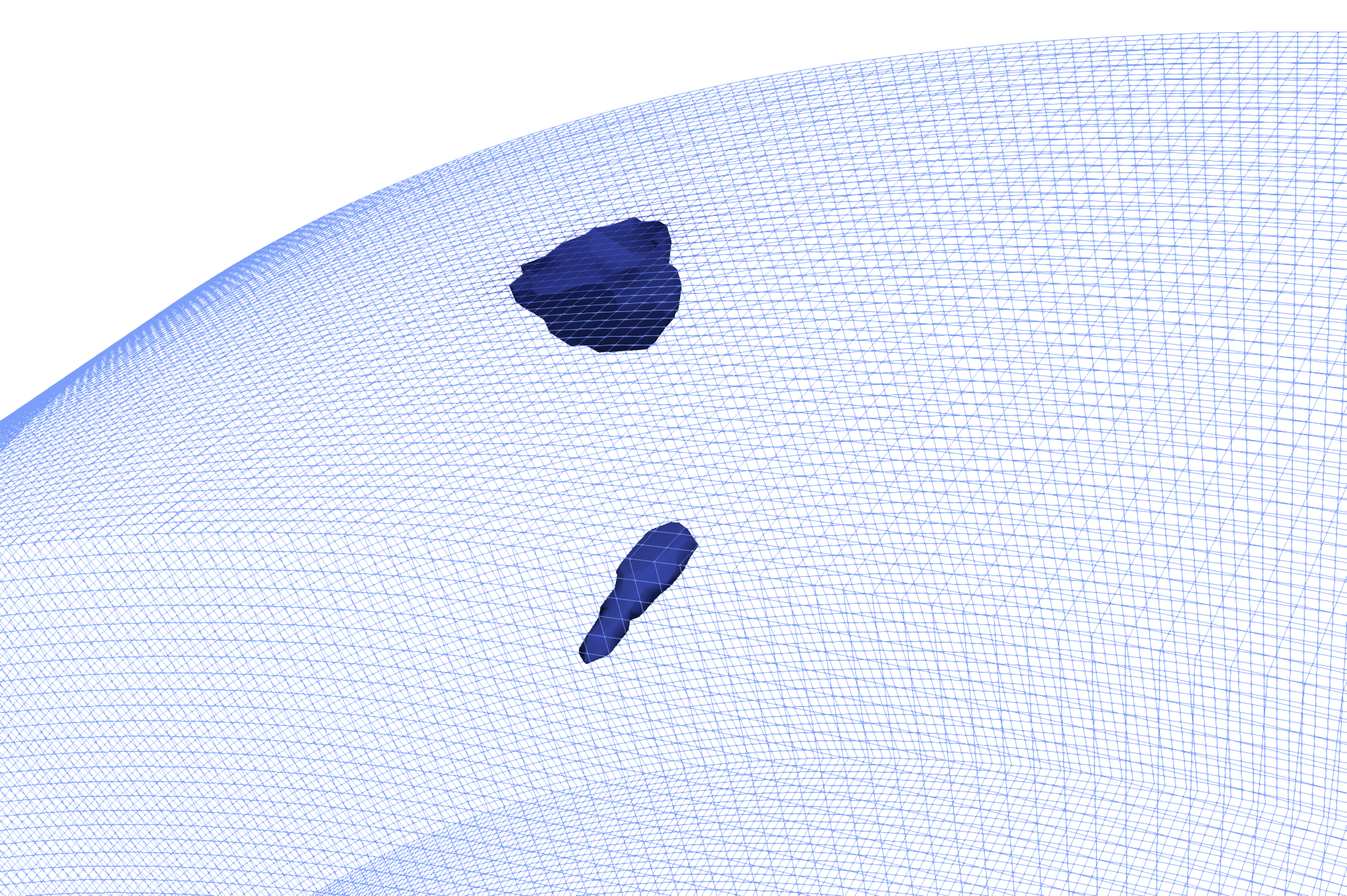}
    &\includegraphics[width=0.18\textwidth]{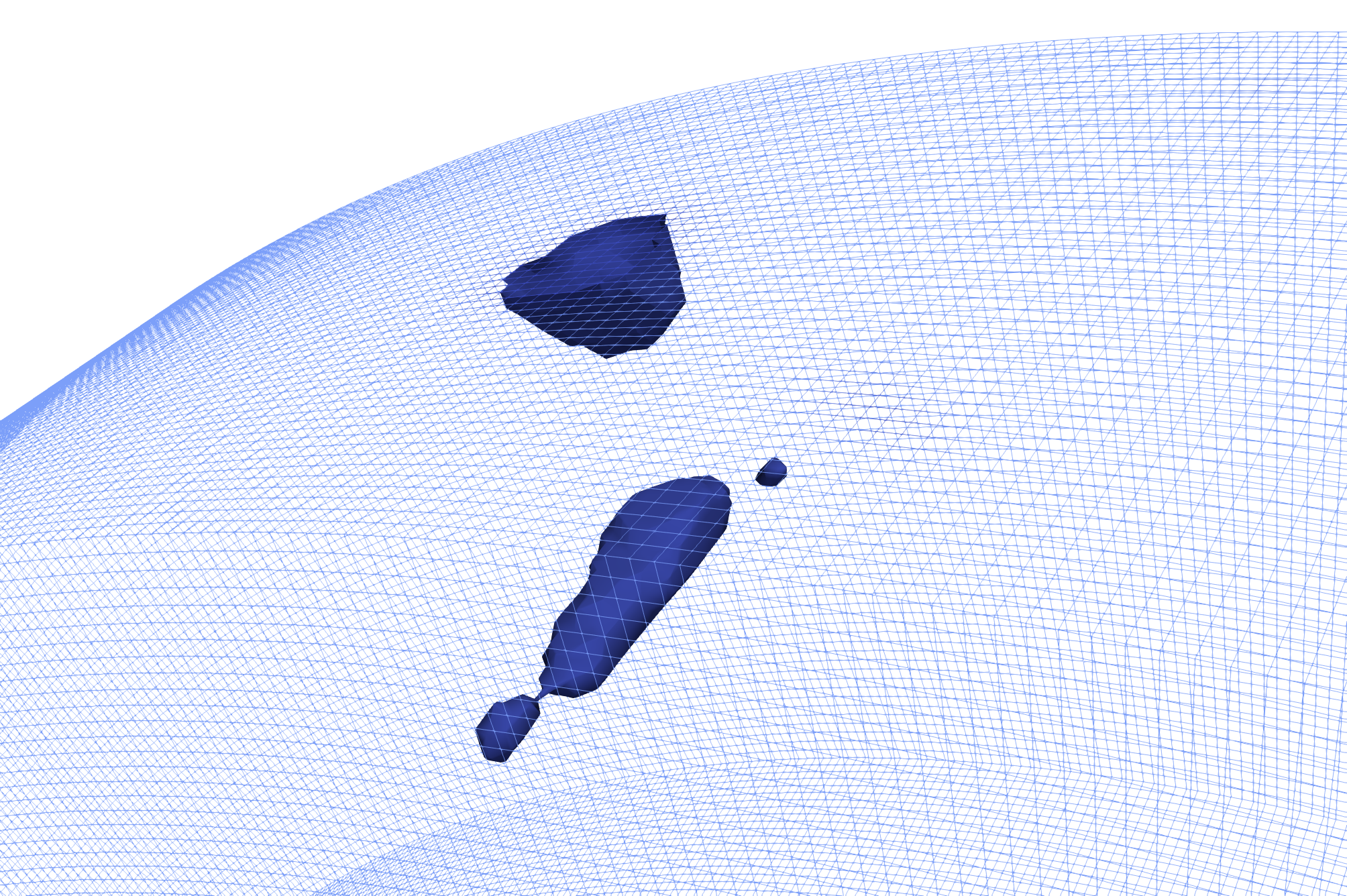}
    &\includegraphics[width=0.18\textwidth]{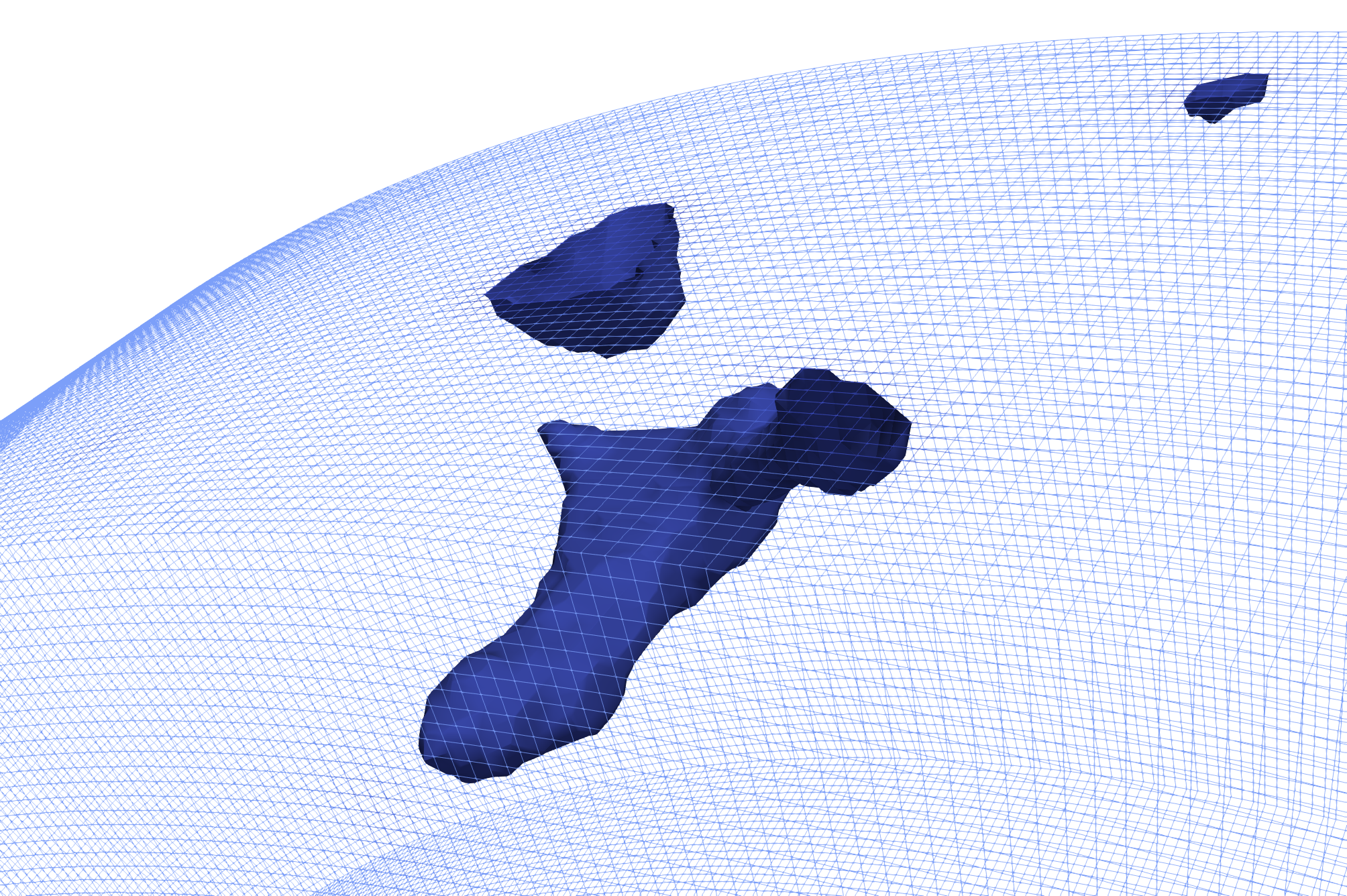}
    &\includegraphics[width=0.18\textwidth]{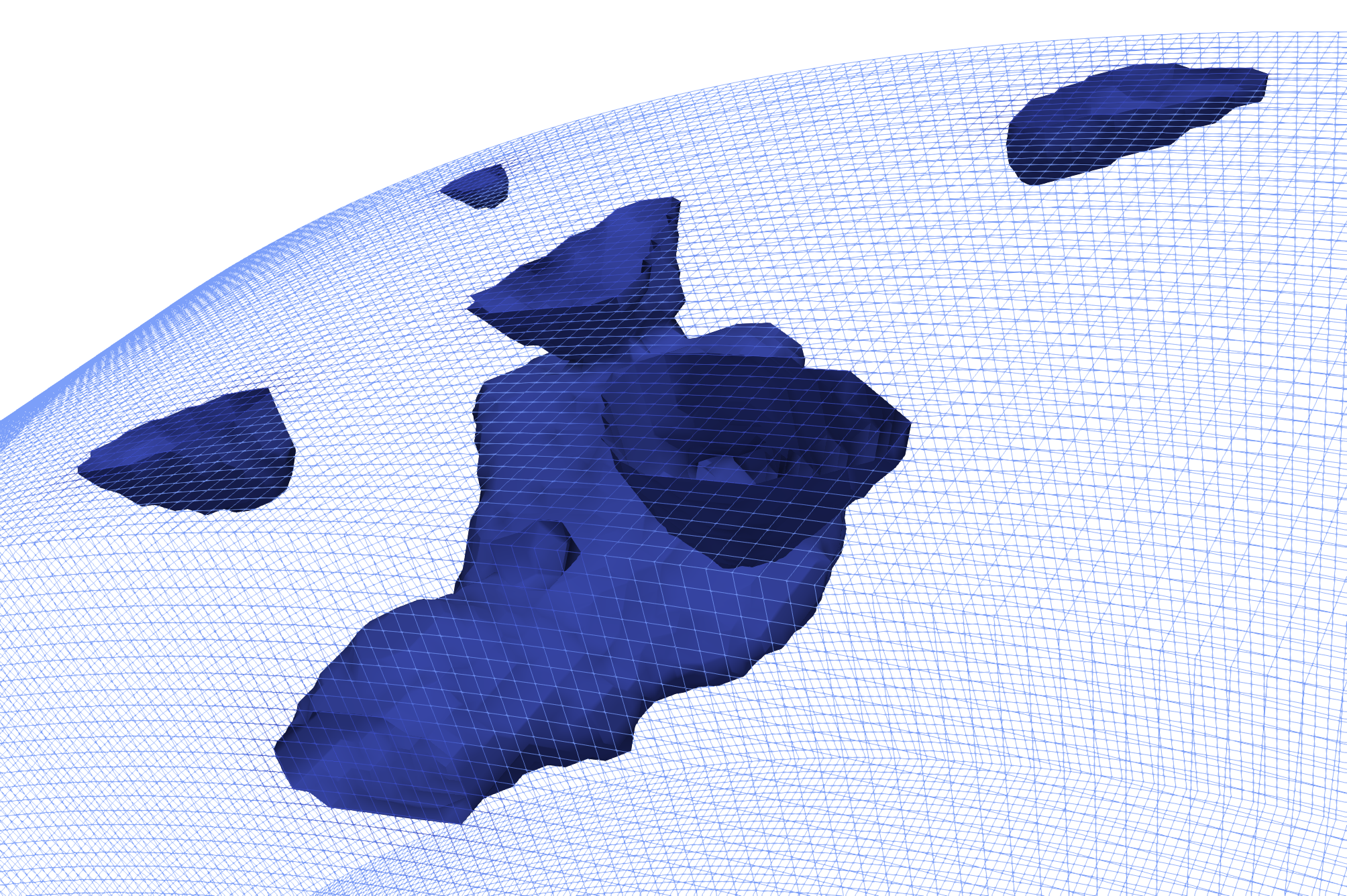}
    \\
  $t=0.65$
  &
  $t=0.6575$
  &
  $t=0.66$
  &
  $t=0.665$
  &
    $t=0.68$
    \\
    \includegraphics[width=0.18\textwidth]{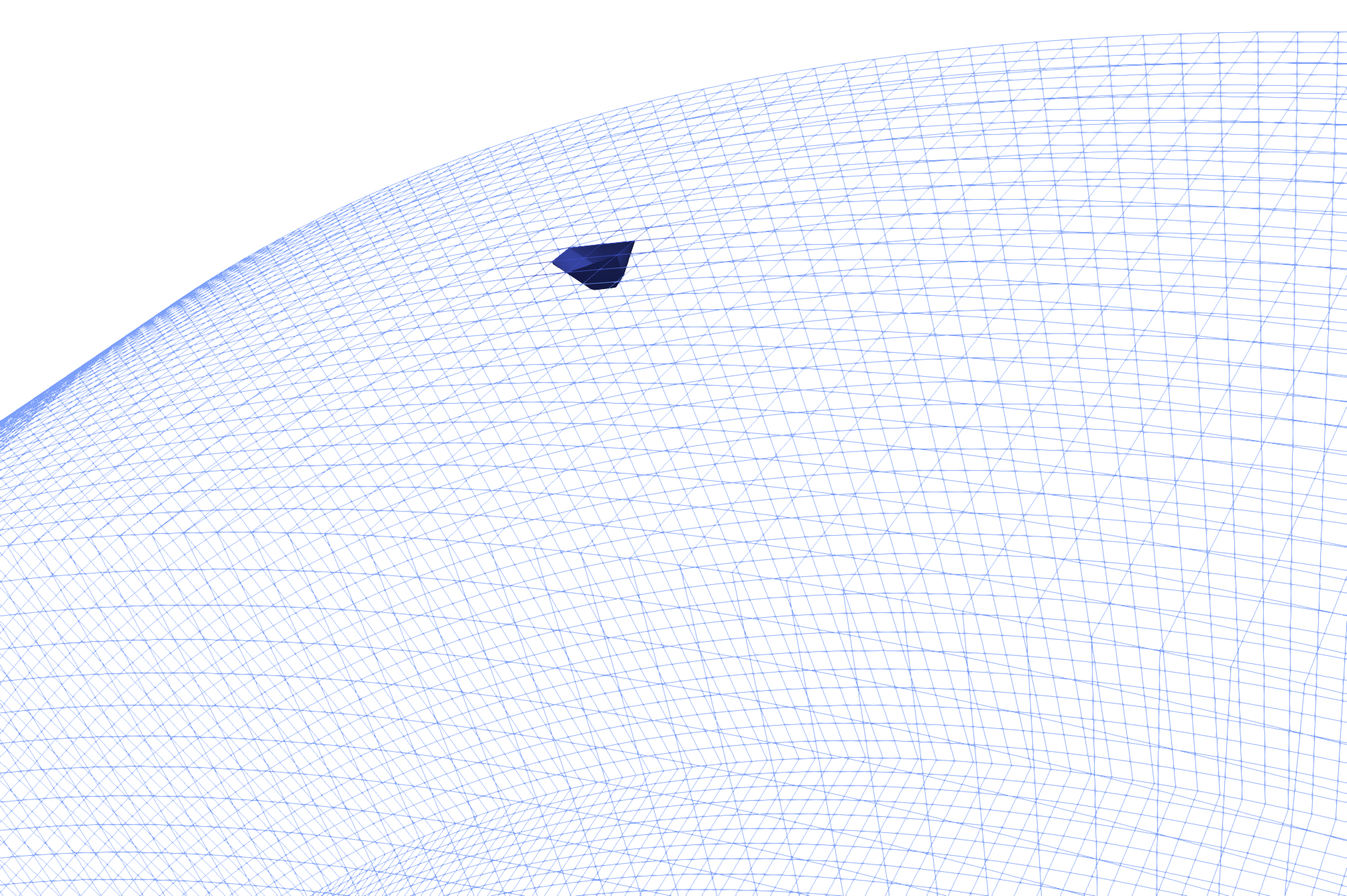}
    &\includegraphics[width=0.18\textwidth]{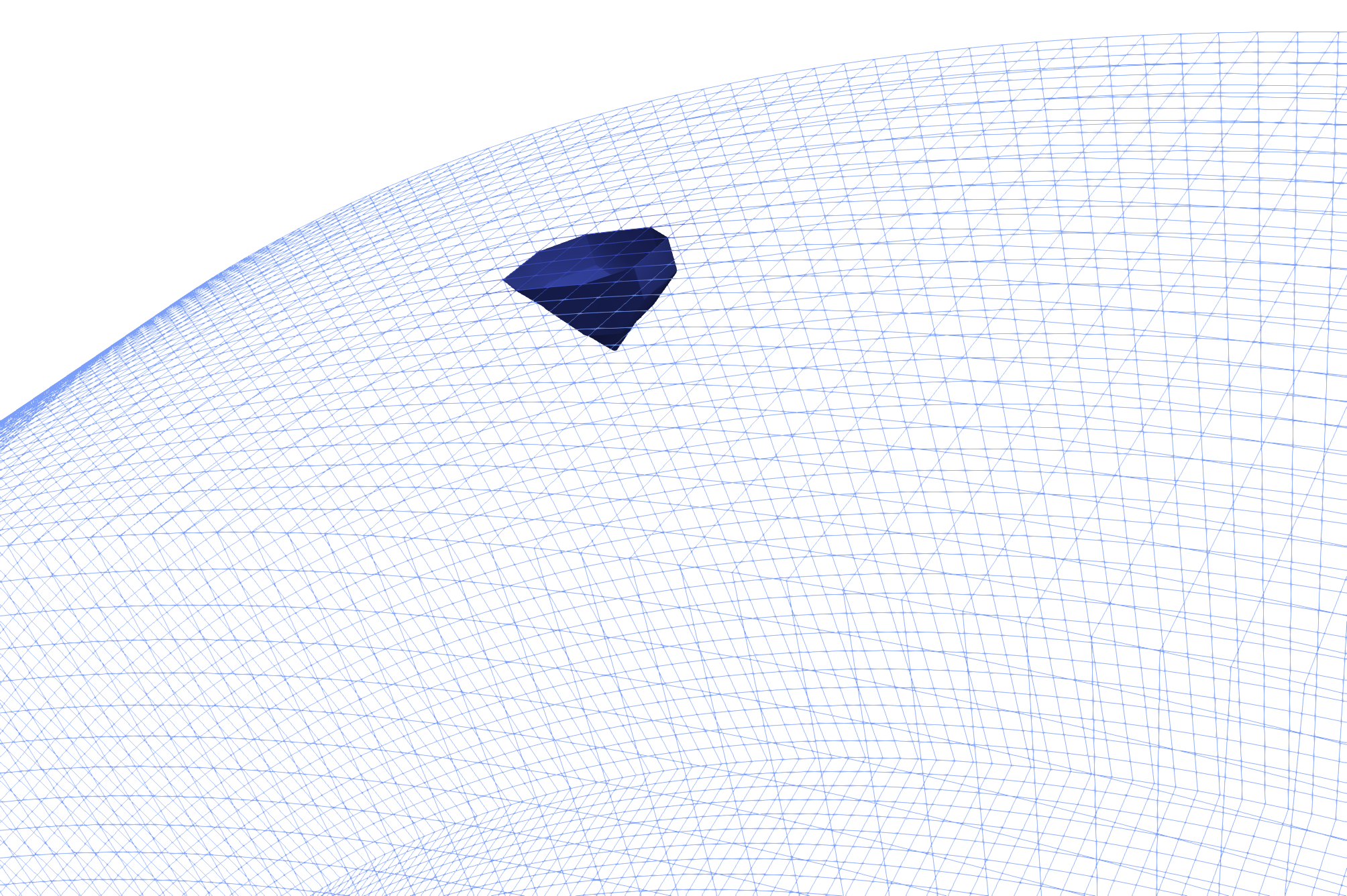}
    &\includegraphics[width=0.18\textwidth]{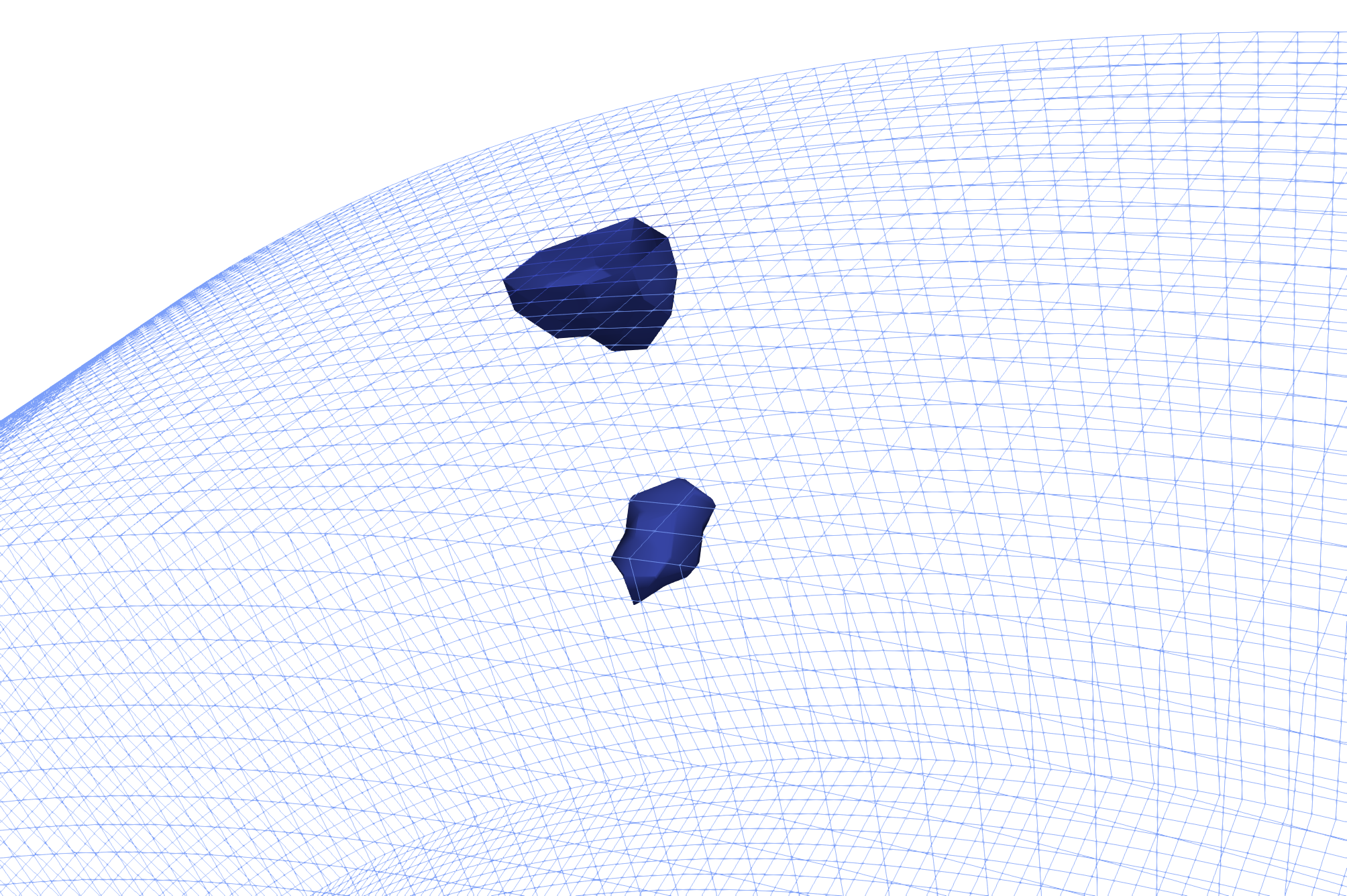}
    &\includegraphics[width=0.18\textwidth]{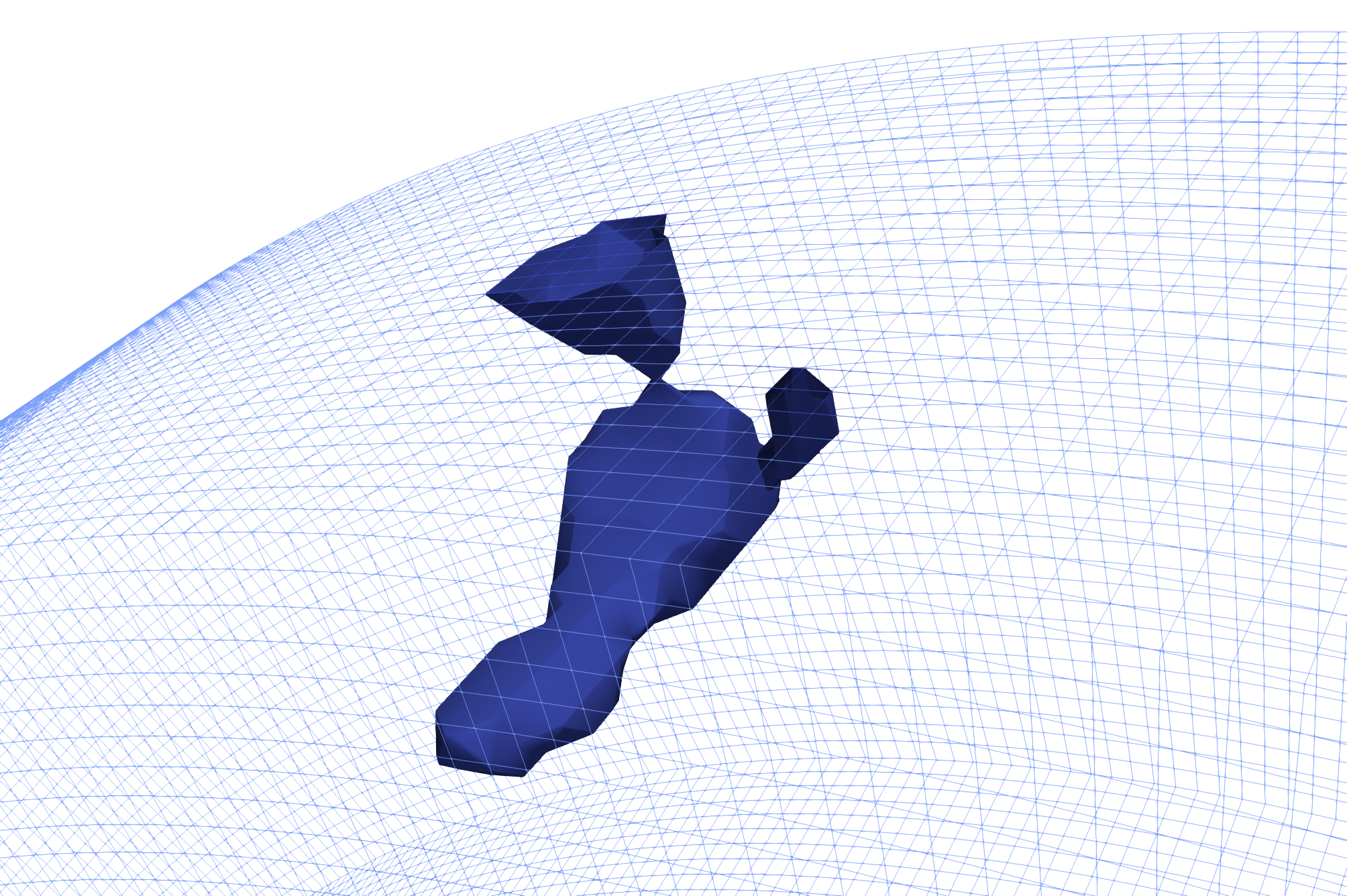}
    &\includegraphics[width=0.18\textwidth]{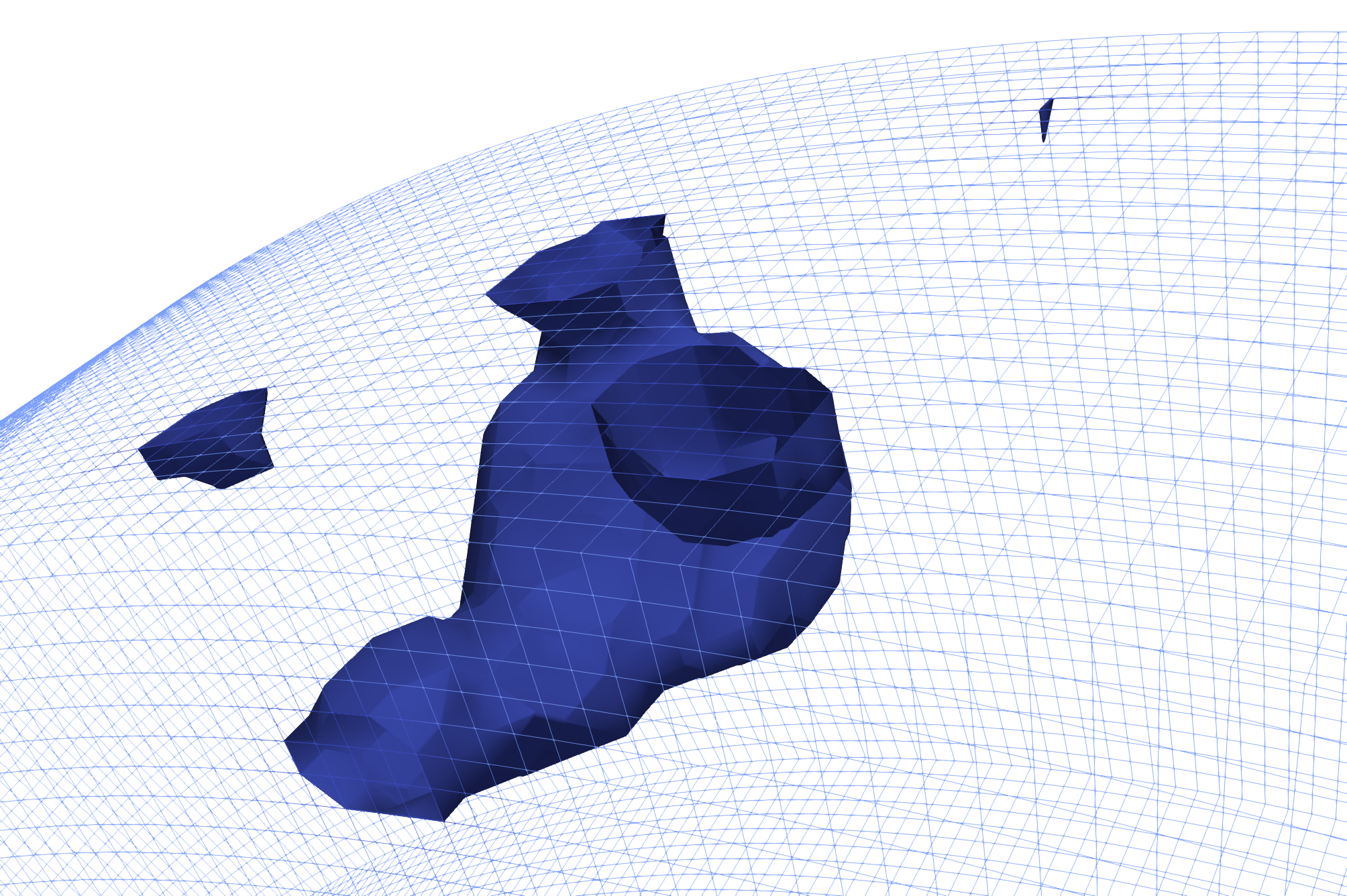}
    \\[2mm]
    \multicolumn{5}{c}{phase-field approximation 
      with  140\,481 DoFs on  131\,072 hexahedra}
\end{tabular}
  \caption{Details of the phase-field evolution at different times with different
mesh resolution.}
\label{fig:3dcracks}
\end{figure}

\gruen{The resulting phase-field approximation
at the final time $t_n= 0.75$ is shown in Fig.~\ref{fig:3dcrack}, and
the evolution of the phase-field approximation at different time steps
is illustrated in Fig.~\ref{fig:3dcracks}.  Comparing the numerical approximation of the phase field on different meshes, we observe that the crack evolution can be described on the coarser level, but more details are included in the finer resolution. On coarser levels, the complexity of the crack pattern cannot be resolved, so that indeed such detailed computations are required to obtain fine fragment spallation.}

Again, we observe wave propagation and crack pattern interaction as a
purely dynamic phenomenon that is only driven by the stress
interaction and does not require any geometric initiation. Our new
first-order method for wave propagation is clearly able to recover
spallation phenomena.

An overview of the software design is reported in \cite{baumgarten2020parallel}, and
the code is available on
\href{https://git.scc.kit.edu/mpp/dgwave/-/tree/fracture1.1}{https://git.scc.kit.edu/mpp/dgwave/-/tree/fracture1.1} including the configuration for the presented examples.

\Clearpage

\section{Conclusion and future perspectives}\label{sec:conclusion}

In this contribution, we aim to establish a new algorithmic
methodology for wave propagation with dynamic fracture.
This methodology allows to reliably predict fracture and fragmentation due
to traveling waves and their superposition as well as the wave
reflections at the fracture interfaces.

\gruen{The major difference to established methods in \cite{Miehe2010,borden2012phase}
is the formulation and approximation of wave propagation as a first-order
hyperbolic system for velocity and stress, using techniques that are
well established in fluid dynamics.  It remains an open question whether it is possible to show the presented fracture dynamic with conforming finite element approximations of the displacements, particularly since DG methods can propagate discontinuities initiated by crack opening. Further investigations are required to see if such a fracture dynamic can also be approximated with a Newmark-type approach.
}

We focus here on the basic formulation for linear elastic materials;
this will be extended to more general material classes and
crack-driving forces in a next step.

In our examples we employ a crack-driving force based on
the {maximum principal stress} $ Y_\text{el} = \displaystyle
\max\Big\{\frac{\sigma_\text{I}}{\sigma_\text{c}} -1,0\Big\}$. The
corresponding state of tension defines the (reversible) phase-field
evolution
    \begin{align*}
      0 \in \tau_\text{r} \dot s +
      Y_\text{el} - M_\text{geom}\partial \gamma_\text{c}(s)
      + \partial \chi_{[0,1]}(s)
      \,,
    \end{align*}
and  transfers to the irreversible phase-field fracture evolution
    \begin{align*}
      0 \in \tau_\text{r} \dot s + \partial \chi_{(-\infty,0)}(\dot s) +
      Y_\text{el} - M_\text{geom}\partial \gamma_\text{c}(s)
      \,.
    \end{align*}
The maximum principal stress driving force is a common fracture
criterion and a straightforward choice for brittle materials.
Likewise, the usual variational formulation of the phase-field driving
force could be employed. This formulation requires the split of the
energy density into a \emph{tensile} and a \emph{compressive} energy
functional $ W_\text{pf}(\vec \varepsilon,s)$, so that the stress
response and phase-field driving force are conjugate and defined by
    \begin{align*}
      \vec \sigma = \partial_{\vec \varepsilon}
      W_\text{pf}(\vec \varepsilon,s)
      \,,\qquad
      Y_\text{el} = -
      \partial_s W_\text{pf}(\vec \varepsilon,s)
      \,.
    \end{align*}
The phase-field fracture formulation can then be extended to problems with finite deformations, see
\cite{hesch2017framework_long,weinberg2016modeling,thomas2017}. However,
it remains an open question how to construct the corresponding discontinuous Galerkin framework.
Another topic of ongoing research is the extension of our scheme to viscoelastic media, as it is analyzed in
\cite{thomas2020approximation} and formulated in a phase-field setting in \cite{bartels2020approximation},
since wave propagation always is dispersive in natural media.

\Paragraph{Acknowledgement}
The authors gratefully acknowledge the support of the Deutsche
Forschungsgemeinschaft (DFG) within the Priority Program~2256 ``Variational Methods for Predicting Complex Phenomena in Engineering Structures and Materials''
in the projects \mbox{WE~2525/15-1} and \mbox{WI~1430/9-1}.


\end{document}